\documentclass[10pt,journal]{IEEEtran}
\usepackage{graphicx, xcolor} 
\usepackage{amsmath, amsfonts, amssymb, amsthm}

\usepackage{xargs}

\newtheorem{theorem}{\textbf{Theorem}}
\newtheorem{lemma}{\textbf{Lemma}}

\usepackage{algorithm,algpseudocode}
\usepackage{algorithmicx}
\usepackage{booktabs}
\usepackage{amsfonts}
\usepackage{comment}
\usepackage[nice]{nicefrac}
\usepackage{subcaption}
\usepackage{xspace}
\usepackage{dsfont}
\usepackage{parskip}

\begin{document}
\title{Quickest Change Detection with Confusing Change}
\author{
    \IEEEauthorblockN{Y.-Z. Janice Chen\IEEEauthorrefmark{1}, Jinhang Zuo\IEEEauthorrefmark{1}, Venugopal V. Veeravalli\IEEEauthorrefmark{3}, Don Towsley\IEEEauthorrefmark{1}}\\
    \IEEEauthorblockA{\IEEEauthorrefmark{1}University of Massachusetts Amherst, \quad\{yuzhenchen, jinhangzuo, towsley\}@cs.umass.edu}\\
    \IEEEauthorblockA{\IEEEauthorrefmark{3}University of Illinois at Urbana-Champaign, \quad vvv@illinois.edu}
    \thanks{This work has been submitted to the IEEE for possible publication. Copyright may be transferred without notice, after which this version may no longer be accessible.}
}

\maketitle

\begin{abstract}
In the problem of quickest change detection (QCD), a change occurs at some unknown time in the distribution of a sequence of independent observations. 
This work studies a QCD problem where the change is either a \textit{bad change}, which we aim to detect, or a \textit{confusing change}, which is not of our interest. 
Our objective is to detect a \textit{bad change} as quickly as possible while avoiding raising a false alarm for pre-change or a \textit{confusing change}. 
We identify a specific set of pre-change, \textit{bad change}, and \textit{confusing change} distributions that pose challenges beyond the capabilities of standard Cumulative Sum (CuSum) procedures.
Proposing novel CuSum-based detection procedures, \texttt{S-CuSum} and \texttt{J-CuSum}, leveraging two CuSum statistics, we offer solutions applicable across all kinds of pre-change, \textit{bad change}, and \textit{confusing change} distributions.
For both \texttt{S-CuSum} and \texttt{J-CuSum}, we provide analytical performance guarantees and validate them by numerical results. 
Furthermore, both procedures are computationally efficient as they only require simple recursive updates. 
\end{abstract}



\section{Introduction}

The quickest change detection (QCD) problem is of fundamental importance in sequential analysis and statistical inference. See, e.g.,~\cite{poor2008quickest, veeravalli2014quickest, tartakovsky2014sequential, xie2021sequential} for books and survey papers on this topic. 
Moreover, the problem of sequential detection of changes or anomalies in stochastic systems arises in a variety of science and engineering domains such as signal processing in sensor networks~\cite{banerjee2015data, sun2022quickest} and cognitive radio~\cite{lai2008quickest}, quality control in manufacturing ~\cite{lai1995sequential, woodall2004using} and power delivery~\cite{banerjee2014power, sun2018backup}, anomaly detection in social network~\cite{ji2017algorithmic} and public health~\cite{liang2023quickest}, and surveillance~\cite{tartakovsky2006novel}, where timely detection of changes or abnormalities is crucial for decision-making or system control.

In the classical formulation of the QCD problem, the objective is to detect \textit{any} change in the distribution. However, in many applications, there can be \textit{confusing changes} that are not of primary interest, and raising an alarm for a \textit{confusing change} may result in wasting human resources on checking the system.  
For instance, in cybersecurity applications, natural fluctuations or diurnal variations in network traffic patterns can be confusing when it comes to identifying potentially malicious activity. In environment or industrial process monitoring, delicate sensors or instruments can experience wear or drift in calibration and, therefore, introduce changes in the data that are not reflective of the underlying environment or process. 
Another example is wireless communication in a complicated environment, where physical obstacles or unrelated nearby communication systems can cause signal interference. 

This study addresses a QCD problem with the aim of swiftly detecting \textit{bad changes} while preventing false alarms for pre-change stages or \textit{confusing changes}. We consider a sequence of observations $X_1, X_2,..., X_t$ generated from a stochastic system. During the pre-change stage (i.e., when the system is in the in-control state), the observations follow distribution $f_0$. At an unknown deterministic time $\nu$, an event occurs and changes the distribution from which the observations arise. The event can manifest as either a confusing change with distribution $f_C$ or a bad change with distribution $f_B$. Addressing the QCD problem with confusing change prompts intriguing questions regarding the design of a detection procedure that will not be triggered by a confusing change. Furthermore, the false alarm metric in this context differs from that in classical QCD problems, necessitating considering both pre-change and confusing change factors.

There has been prior research exploring extensions of the classical QCD framework to encompass formulations with more intricate assumptions regarding distributions before or after an event occurs, including composite pre-change distribution~\cite{mei2003asymptotically, mei2006sequential}, composite post-change distribution~\cite{rovatsos2017statistical}, post-change distribution isolation~\cite{warner2022sequential, zhao2022adaptive}, and transient change~\cite{zou2018quickest}. 
Our QCD with confusing change problems differ from those problems as they still aim to raise an alarm as soon as \textit{any} change happens, while we avoid raising an alarm if the change is a confusing change. 
The most relevant extension for our study is the formulation of nuisance change, as studied by Lau and Tay~\cite{lau2019quickest}. In a similar vein,~\cite{lau2019quickest} considers two types of changes and aims to detect only one of them. 
While~\cite{lau2019quickest} presents a more generalized model for change points, allowing one type of change to occur after the other, they focus solely on a restricted set of post-change distributions. In contrast to~\cite{lau2019quickest}, we focus on scenarios where one type of change occurs and devise solutions applicable to all post-change distributions. Therefore, we view our work as complementary to Lau and Tay's~\cite{lau2019quickest}. 
From an algorithm design perspective, our proposed procedures differ from prior procedures that involve two CuSum Statistics~\cite{dragalin1997design, zhao2005dual}, and our procedures are more computationally efficient than the window-limited generalized likelihood ratio test-based procedure proposed in~\cite{lau2019quickest}.

Our contributions are as follows: In Section~\ref{sec:formulation}, we formulate a novel quickest change detection problem where the change can be either a bad change or a confusing change, and we propose a false alarm metric that takes not only the run length to false alarm for pre-change but also the run length to false alarm for confusing change into account. In Section~\ref{sec:scenario}, we investigate the problem with possible scenarios with different combinations of pre-change, bad change, and confusing change distributions, and we identify a scenario in which procedures depending on a single CuSum statistic will incur short run lengths to false alarm. In Section~\ref{sec:algos}, we propose two novel procedures, \texttt{S-CuSum} and \texttt{J-CuSum}, that incorporate two CuSum statistics and work under all kinds of distributions. In Section~\ref{sec:analysis}, we first provide a universal lower bound of detection delay for any procedure that fulfills the false alarm requirement, and then we provide the false alarm upper bounds and detection delay lower bounds of \texttt{S-CuSum} and \texttt{J-CuSum}. In Section~\ref{sec:simulations}, our simulation results under all possible scenarios corroborate our theoretical guarantees for \texttt{S-CuSum} and \texttt{J-CuSum}. Finally, in Section~\ref{sec:conclusion}, we conclude this paper and discuss future directions.


\section{Problem Formulation \& Preliminary}\label{sec:formulation}

Let $\{X_t: t \in \mathbb{N}_{+}\}$ be a sequence of random variables whose values are observed sequentially, and let $\{\mathcal{F}_t: t\in \mathbb{N}_{+}\}$ be the filtration generated by the sequence, i.e., $\mathcal{F}_t := \sigma(X_i:1\leq t\leq t)$, where $\sigma(\cdot)$ denotes $\sigma$-algebra. Let $f_0(X)$ be the density of the pre-change distribution, let $f_C(X)$ be the density of the confusing change distribution, and let $f_B(X)$ be the density of the bad change distribution. We assume that $f_0$, $f_C$, and $f_B$ are known, different, and measure over a common measurable space. We denote the unknown but deterministic change-point by $\nu$. 

To be more specific, we denote by $\mathbb{P}_{\nu, f_*}$, $*\in \{C, B\}$, the underlying probability measure, and by $\mathbb{E}_{\nu, f_*}$ the corresponding expectation, when the change-point is $\nu$ and the post-change distribution is $f_*$. Moreover, we denote by $\mathbb{P}_{\infty}$the underlying probability measure, and by $\mathbb{E}_{\infty}$ the corresponding expectation, when the change never occurs, i.e., $\mathbb{P}_{\infty} = \mathbb{P}_{f_0}$, $\mathbb{E}_{\infty} = \mathbb{E}_{f_0}$. Finally, let $T$ denote a stopping time, i.e., the time at which we stop taking observations and declare that a bad change has occurred.

\textbf{False Alarm Measure.} An alarm is considered as a false alarm if (1) the alarm is raised in the pre-change stage, i.e., $T < \nu$, or (2) the alarm is raised for a confusing change, i.e., $T \geq \nu$ and $X_{T} \sim f_{C}(X)$. Therefore, for any stopping time $T$, we measure the false alarm performance in terms of its mean time to false alarm in both cases. Specifically, we denote by $\mathcal{C}_{\gamma}$ the subfamily of stopping times for which
the worst-case average run length (WARL) to false alarm
is at least $\gamma$ in both cases, i.e.,
\begin{align}\label{eq:stopping-time-set}
    \mathcal{C}_{\gamma} := \{T: \mathbb{E}_{\infty}[T] \geq \gamma, \inf_{\nu\geq 1}\mathbb{E}_{\nu, f_C}[T] \geq \gamma\}.
\end{align}

\textbf{Delay Measure.} We use worst-case measures for delay.
By Pollak’s criterion~\cite{pollak1985optimal}, for any $T\in \mathcal{C}_{\gamma}$, 
\begin{align}
    \text{WADD}_{f_B}(T) := \sup\limits_{\nu \geq 1}\mathbb{E}_{\nu, f_B}[T-\nu|T \geq \nu].
\end{align}

\textbf{Optimization Problem.} The optimization problem we have in the work is to find a stopping rule for a bad change that belongs to $\mathcal{C}_{\gamma}$ for every $\gamma \geq 1$ and approximates 
\begin{align}
    &\inf\limits_{T\in \mathcal{C}_{\gamma}}\text{WADD}_{f_B}(T) 
\end{align}
as $\gamma \rightarrow \infty$.

\textbf{Methodology.} In this work, we propose Cumulative Sum (CuSum)-based procedures for the QCD with confusing change problems. Hence, we briefly review CuSum statistics and its recursion here. We let the cumulative summation of the log-likelihood ratio of two densities, e.g., $f_B$ and $f_0$, be:
\begin{align}\label{eq:CuSumfBf0}
    \text{CuSum}_{f_B, f_0}[t] := \max\limits_{k:1\leq k\leq t} \sum_{i=k}^t \log\left(\frac{f_B(X_i)}{f_0(X_i)}\right).
\end{align}
Such a CuSum statistic can be updated recursively by:
\begin{align}
    \text{CuSum}_{f_B, f_0}[t]\equiv \left(\text{CuSum}_{f_B, f_0}[t-1] + \log\left(\frac{f_B(X_t)}{f_0(X_t)}\right)\right)^{+}. 
\end{align}
Please refer to, e.g.,~\cite{poor2008quickest, veeravalli2014quickest, tartakovsky2014sequential, xie2021sequential} for more discussion on this type of procedure.


\begin{table}[t]
\centering
\caption{Possible Drifts of $\text{CuSum}_{f_B, f_0}$ and $\text{CuSum}_{f_B, f_C}$ Statistics under Distributions $f_0$, $f_C$, and $f_B$}
\begin{tabular}{ll}
\toprule
 &  $\text{CuSum}_{f_B, f_0}$   \\
\midrule
Under $f_0$ & $\mathbb{E}_{f_0}[\log(\nicefrac{f_B(X)}{f_0(X)})]$\\
& $=-D_{KL}(f_0||f_B) < 0$\\
\midrule
Under $f_C$ & $E_{f_C} [\log(\nicefrac{f_B (X)}{f_0(X))}]$\\
&$= -D_{KL}(f_C||f_B) +D_{KL}(f_C||f_0) \,\,{\pmb \gtrless }\,\, 0$  \\
\midrule
Under $f_B$ & $\mathbb{E}_{f_B}[\log(\nicefrac{f_B(X)}{f_0(X)})]$\\
&$=D_{KL}(f_B||f_0) > 0$  \\
\bottomrule
\toprule
 &  $\text{CuSum}_{f_B, f_C}$  \\
\midrule
Under $f_0$ & $\mathbb{E}_{f_0}[\log(\nicefrac{f_B(X)}{f_C(X)})]$\\
&$=-D_{KL}(f_0||f_B) + D_{KL}(f_0||f_C) \,\,{\pmb \gtrless}\,\, 0$\\
\midrule
Under $f_C$ & $\mathbb{E}_{f_C}[\log(\nicefrac{f_B(X)}{f_C(X)})]$\\
&$=-D_{KL}(f_C||f_B) < 0$\\
\midrule
Under $f_B$ & $\mathbb{E}_{f_B}[\log(\nicefrac{f_B(X)}{f_C(X)})]$\\
&$=D_{KL}(f_B||f_C) > 0$\\
\bottomrule
\end{tabular} 
\label{tbl:drifts-half2}
\end{table}

\section{Characterization of Challenging Scenario}\label{sec:scenario}

In this section, we categorize all instances of the QCD with confusing change problems into three scenarios and identify one scenario as more challenging than the other two. Specifically, we characterize the scenarios by the values of KL divergences $D_{KL}(f_0||f_C)$, $D_{KL}(f_0||f_B)$, $D_{KL}(f_C||f_0)$, $D_{KL}(f_C||f_B)$, $D_{KL}(f_B||f_0)$, and  $D_{KL}(f_B||f_C)$ as they govern the \textit{drifts} of CuSum statistics.

The drifts of CuSum statistics are determined by the signs of the expectations of the log-likelihood ratios. For example, when a data sample comes from distribution $f_B$, i.e., $X\sim f_B$, the expectation of log-likelihood ratio $\log(f_B(X)/f_0(X))$ corresponds to KL divergence $D_{KL}(f_B||f_0)$, as
\begin{align}
    D_{KL}(f_B||f_0) := \mathbb{E}_{f_B}\left[\log\left(\frac{f_B(X)}{f_0(X)}\right)\right]>0,
\end{align}
where by definition $D_{KL}(f_B||f_0) > 0$ as long as $f_B \neq f_0$, in which case, the value of $\text{CuSum}_{f_B, f_0}$ is generally increasing, i.e., $\text{CuSum}_{f_B, f_0}$ has a positive drift. Following the same reasoning, under $f_0$, $\text{CuSum}_{f_B, f_0}$ has a negative drift because 
\begin{align}
    \mathbb{E}_{f_0}\left[\log\left(\frac{f_B(X)}{f_0(X)}\right)\right] = -D_{KL}(f_B||f_0) < 0.
\end{align}
Interestingly, the drift of $\text{CuSum}_{f_B, f_0}$ under $f_C$, indicated by
\begin{align}
    &\mathbb{E}_{f_C}\left[\log\left(\frac{f_B(X)}{f_0(X)}\right)\right] \notag\\
    &= \mathbb{E}_{f_C}\left[\log\left(\frac{f_B(X)}{f_0(X)}\frac{f_C(X)}{f_C(X)}\right)\right]\\
    &= \mathbb{E}_{f_C}\left[\log\left(\frac{f_B(X)}{f_C(X)}\right)\right] + \mathbb{E}_{f_C}\left[\log\left(\frac{f_C(X)}{f_0(X)}\right)\right] \\
    &= - D_{KL}(f_C||f_B) + D_{KL}(f_C||f_0),
\end{align}
can be positive, zero, or negative depending on the values of $-D_{KL}(f_C||f_B)$ and $D_{KL}(f_C||f_0)$.

As there is an extra distribution, $f_C$, besides the typical distributions, $f_0, f_B$, in this QCD with confusing change problem, it is natural to consider an extra CuSum statistic, $\text{CuSum}_{f_B, f_C}$ (defined in the same way as $\text{CuSum}_{f_B, f_0}$ in \eqref{eq:CuSumfBf0}), besides the typical CuSum statistic, $\text{CuSum}_{f_B, f_0}$, often considered in quickest change detection problem. (Further justification for considering $\text{CuSum}_{f_B, f_C}$ will be provided in Section~\ref{sec:algos}.) Following the same reasoning for $\text{CuSum}_{f_B, f_0}$, we find that $\text{CuSum}_{f_B, f_C}$ has a positive drift under $f_B$, has a negative drift under $f_C$, and has either a positive, zero, or negative drift under $f_0$ depending on the values of $-D_{KL}(f_0||f_B)$ and $D_{KL}(f_0||f_C)$. We list the possible drifts of the two CuSum statistics under the three possible distributions, $f_0, f_B, f_C$, in Table~\ref{tbl:drifts-half2}.

\begin{figure}[t]
    \centering
    \hfill
    \begin{minipage}[b]{0.48\linewidth}
        \includegraphics[width=\columnwidth]{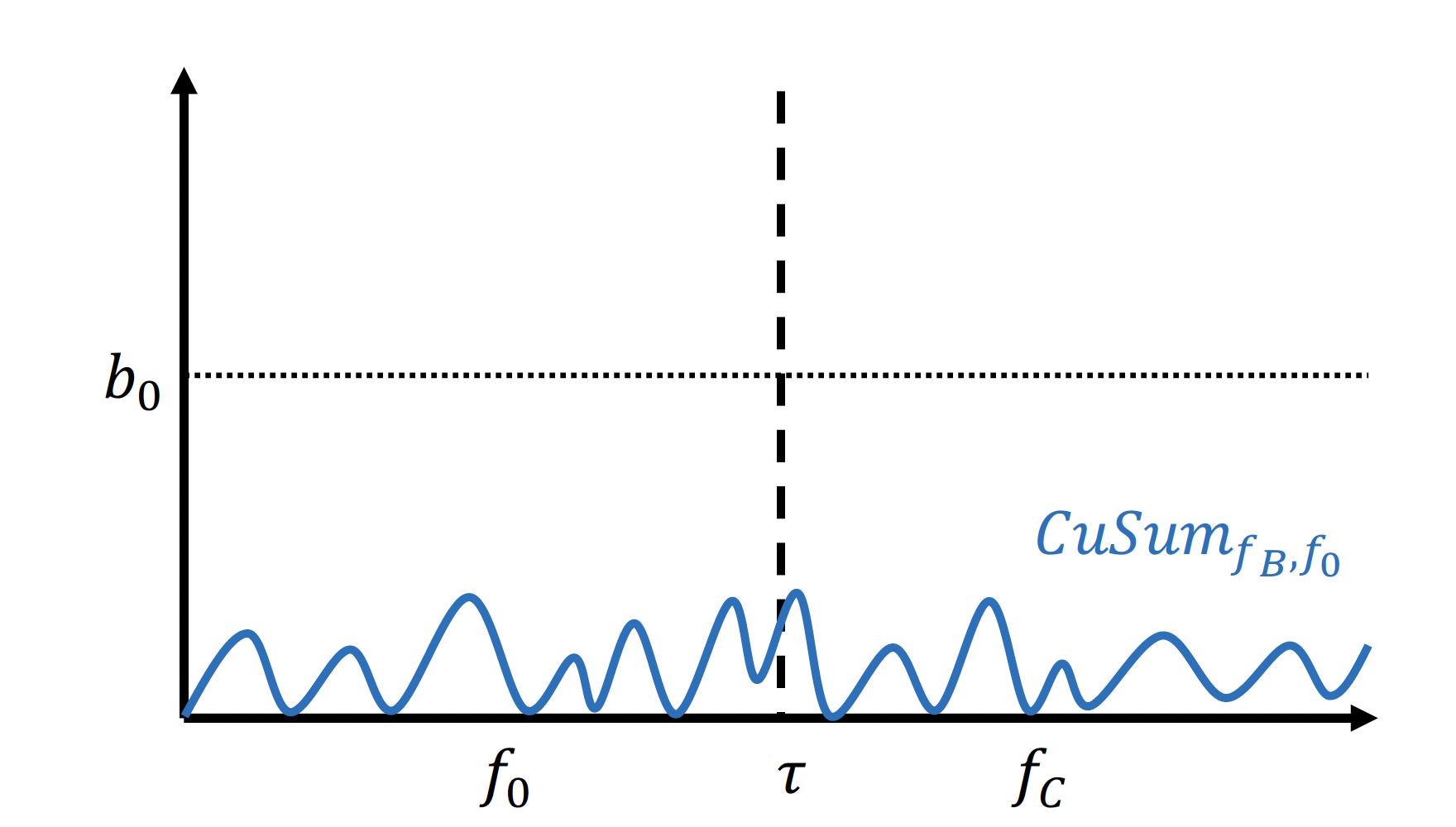}
        \subcaption{Confusing Change}
        \label{fig:trivial-1-confusing}
    \end{minipage}
    \hfill
    \begin{minipage}[b]{0.48\linewidth}
        \includegraphics[width=\columnwidth]{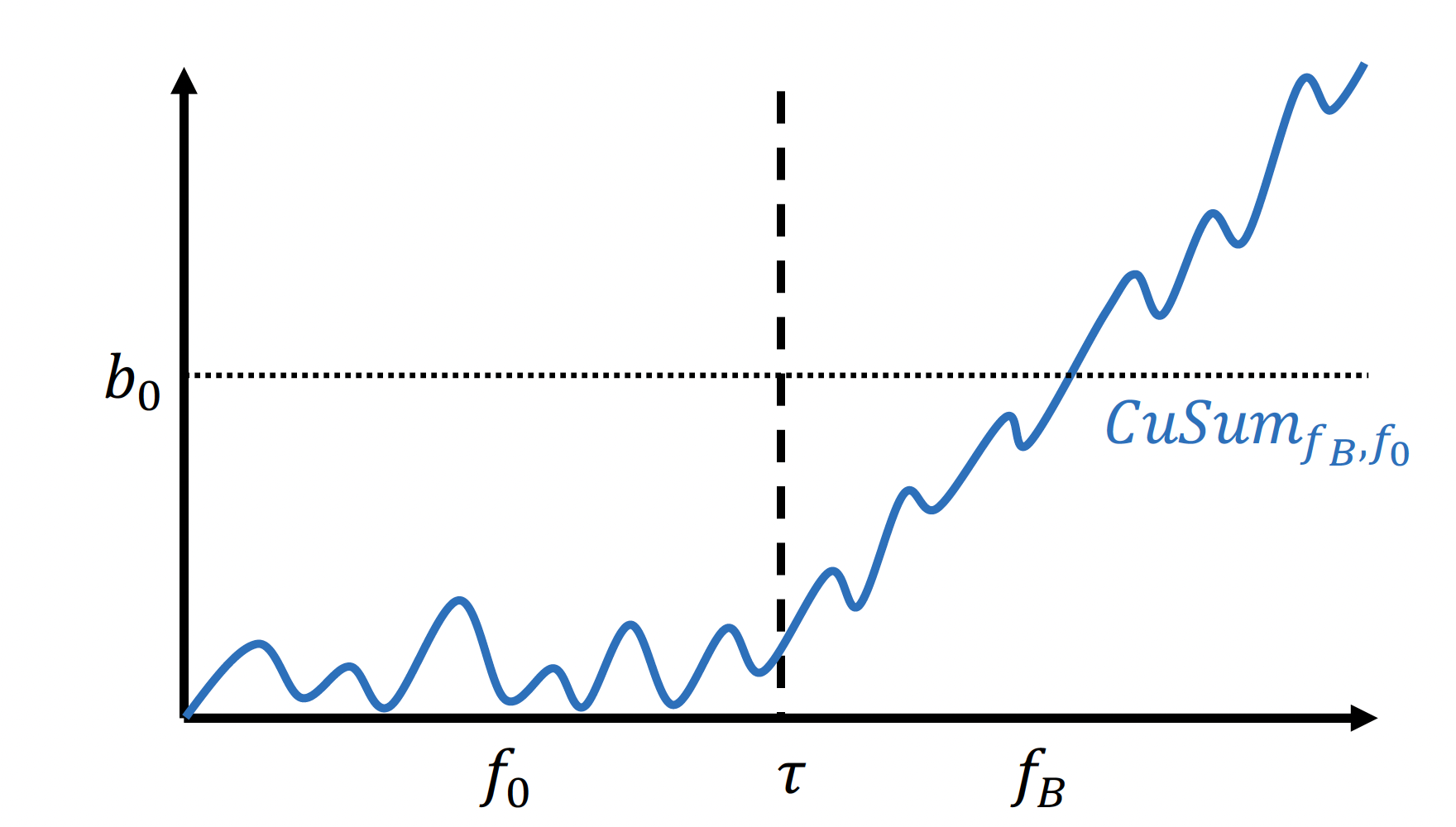}
        \subcaption{Bad Change}
        \label{fig:trivial-1-bad}
    \end{minipage}
    \hfill
    \caption{Illustration of $\text{CuSum}_{f_B, f_0}$ Behaviors in \textbf{Scenario 1} - applying standard $\text{CuSum}_{f_B, f_0}$ procedure suffices to detect a bad change quickly while avoiding raising false alarm for a confusing change.}
    \label{fig:scenario-1}
\end{figure}

\begin{figure}[t]
    \centering
    \hfill
    \begin{minipage}[b]{0.48\linewidth}
        \includegraphics[width=\columnwidth]{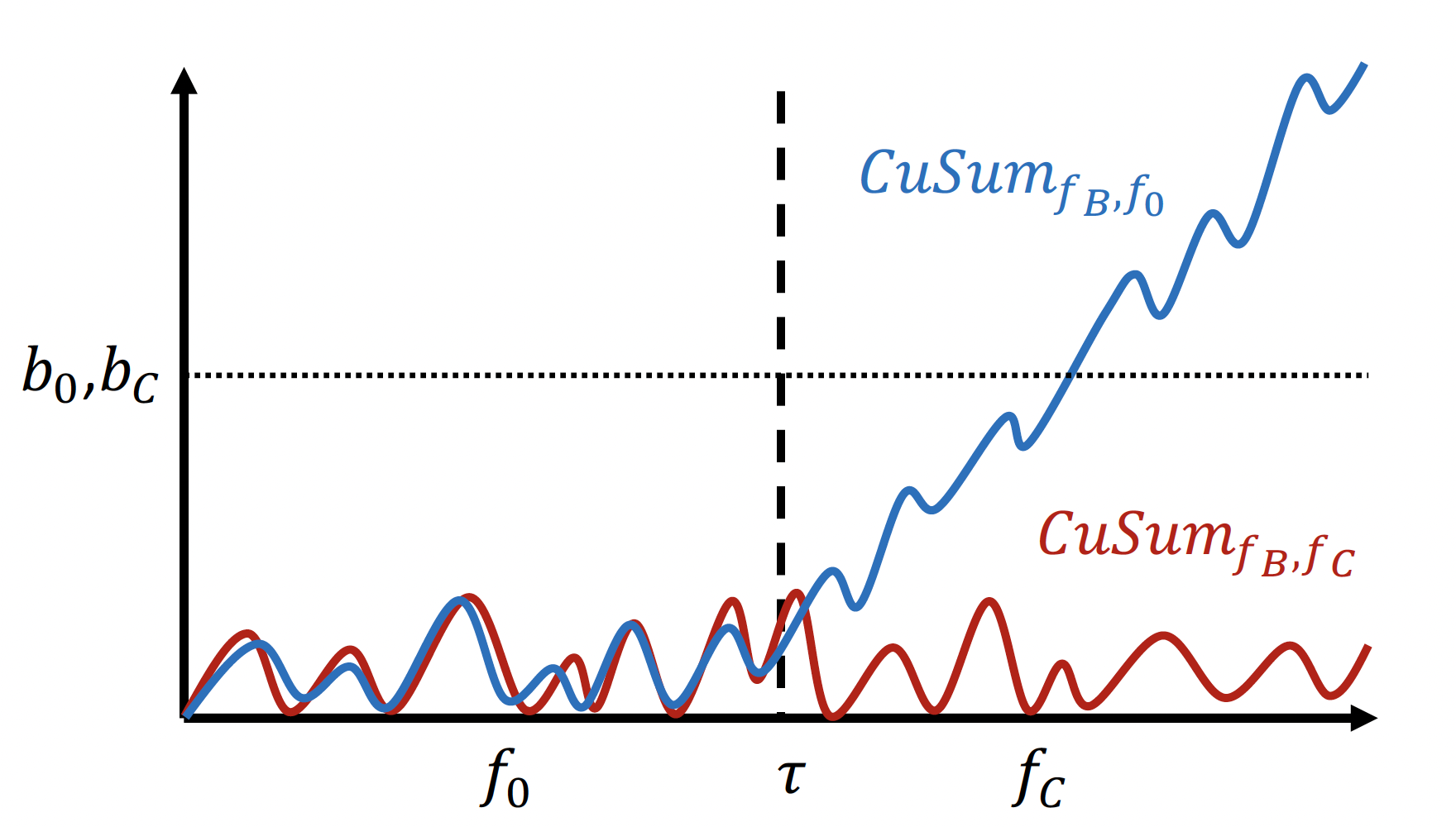}
        \subcaption{Confusing Change}
        \label{fig:trivial-2-confusing}
    \end{minipage}
    \hfill
    \begin{minipage}[b]{0.48\linewidth}
        \includegraphics[width=\columnwidth]{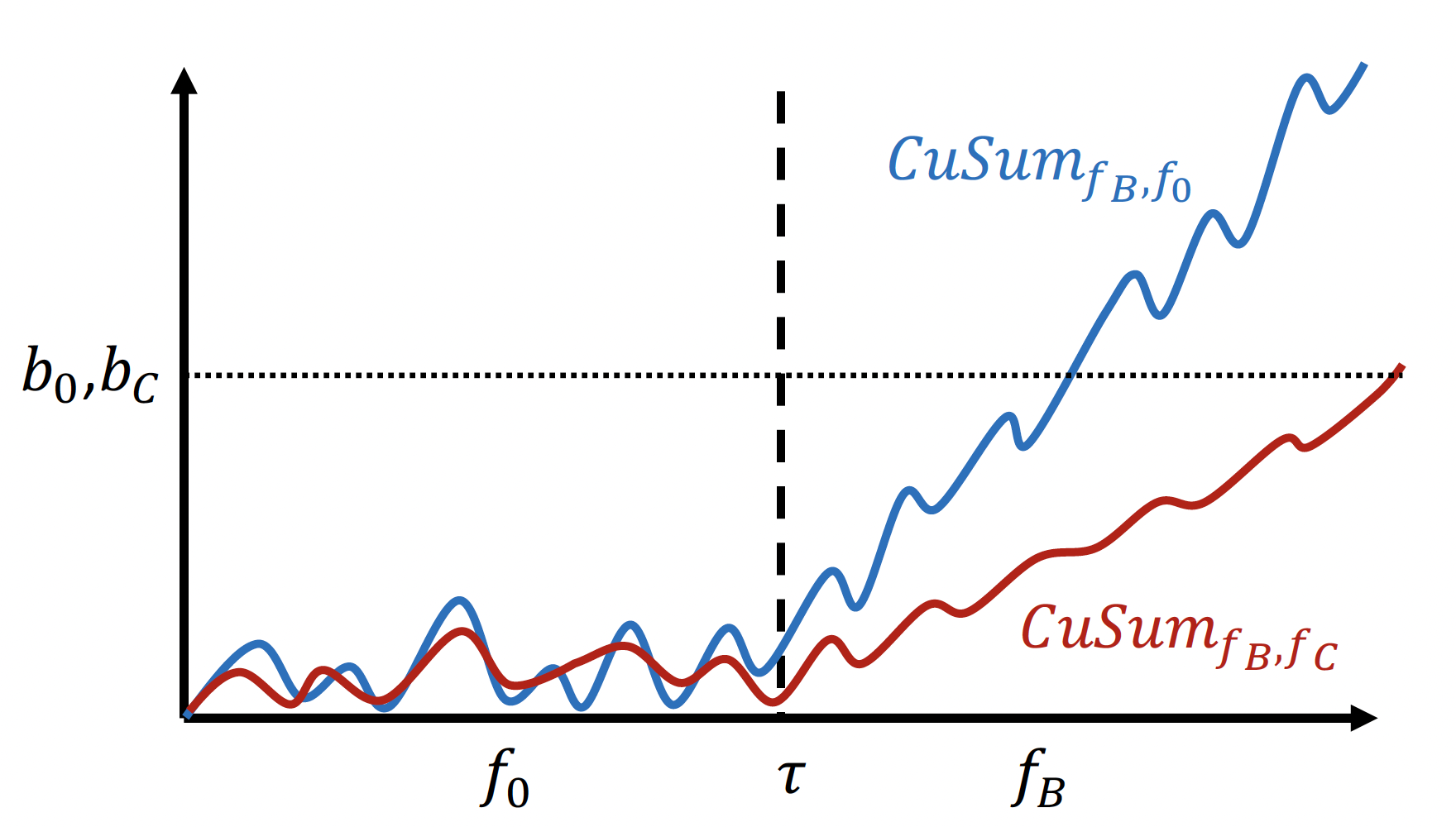}
        \subcaption{Bad Change}
        \label{fig:trivial-2-bad}
    \end{minipage}
    \hfill
    \caption{Illustration of $\text{CuSum}_{f_B, f_0}$ (Blue Lines) and $\text{CuSum}_{f_B, f_C}$ (Red Lines) Behaviors in \textbf{Scenario 2} - applying standard $\text{CuSum}_{f_B, f_C}$ procedure suffices to detect a bad change quickly while avoiding raising false alarm for a confusing change.}
    \label{fig:scenario-2}
\end{figure}

According to the relations between the KL divergences for $f_0$, $f_C$, and $f_B$, which leads to different drifts in the CuSum statistics, we categorize instances of the QCD with confusing change problem into three scenarios: 
\begin{enumerate}
    \item[\textbf{1}.] when $\text{CuSum}_{f_B, f_0}$ has a zero or negative drift under $f_C$, i.e., $-D_{KL}(f_C||f_B) +D_{KL}(f_C||f_0) \leq 0$;
    \item[\textbf{2}.] when $\text{CuSum}_{f_B, f_0}$ has a positive drift under $f_C$ and $\text{CuSum}_{f_B, f_C}$ has a zero or negative drift under $f_0$, i.e., $-D_{KL}(f_C||f_B) +D_{KL}(f_C||f_0) > 0$ and $-D_{KL}(f_0||f_B) + D_{KL}(f_0||f_C) \leq 0$;
    \item[\textbf{3}.] when $\text{CuSum}_{f_B, f_0}$ has a positive drift under $f_C$ and $\text{CuSum}_{f_B, f_C}$ also has a positive drift under $f_0$, i.e., $-D_{KL}(f_C||f_B) +D_{KL}(f_C||f_0) > 0$ and $-D_{KL}(f_0||f_B) + D_{KL}(f_0||f_C) > 0$.
\end{enumerate}

Scenarios 1 and 2 are trivial since simply applying a standard CuSum procedure suffices to quickly detect a bad change while avoiding raising a false alarm for a confusing change. In Figure~\ref{fig:scenario-1}, we illustrate the behavior of $\text{CuSum}_{f_B, f_0}$ in Scenario 1. 
As illustrated in Figure~\ref{fig:trivial-1-confusing}, since in Scenario 1, $\text{CuSum}_{f_B, f_0}$ has a non-positive drift under confusing change $f_C$ (as under pre-change $f_0$), $\text{CuSum}_{f_B, f_0}$ generally remains small after a confusing change occurs and will rarely trigger a false alarm. Figure~\ref{fig:trivial-1-bad} illustrates that $\text{CuSum}_{f_B, f_0}$ generally increases after a bad change occurs as it has a positive drift and quickly passes the alarm threshold. In Figure~\ref{fig:scenario-2}, we find that in Scenario 2 $\text{CuSum}_{f_B, f_0}$ fails to distinguish a bad change from confusing change because it has a positive drift under confusing change $f_C$ (as under bad change $f_B$). Fortunately, applying a standard $\text{CuSum}_{f_B, f_C}$ procedure suffices for our goal as it has non-positive drifts under both the pre-change and confusing change distributions and only has positive drift under the bad change distribution.

\begin{figure}[t]
    \centering
    \hfill
    \begin{minipage}[b]{0.48\linewidth}
        \includegraphics[width=\columnwidth]{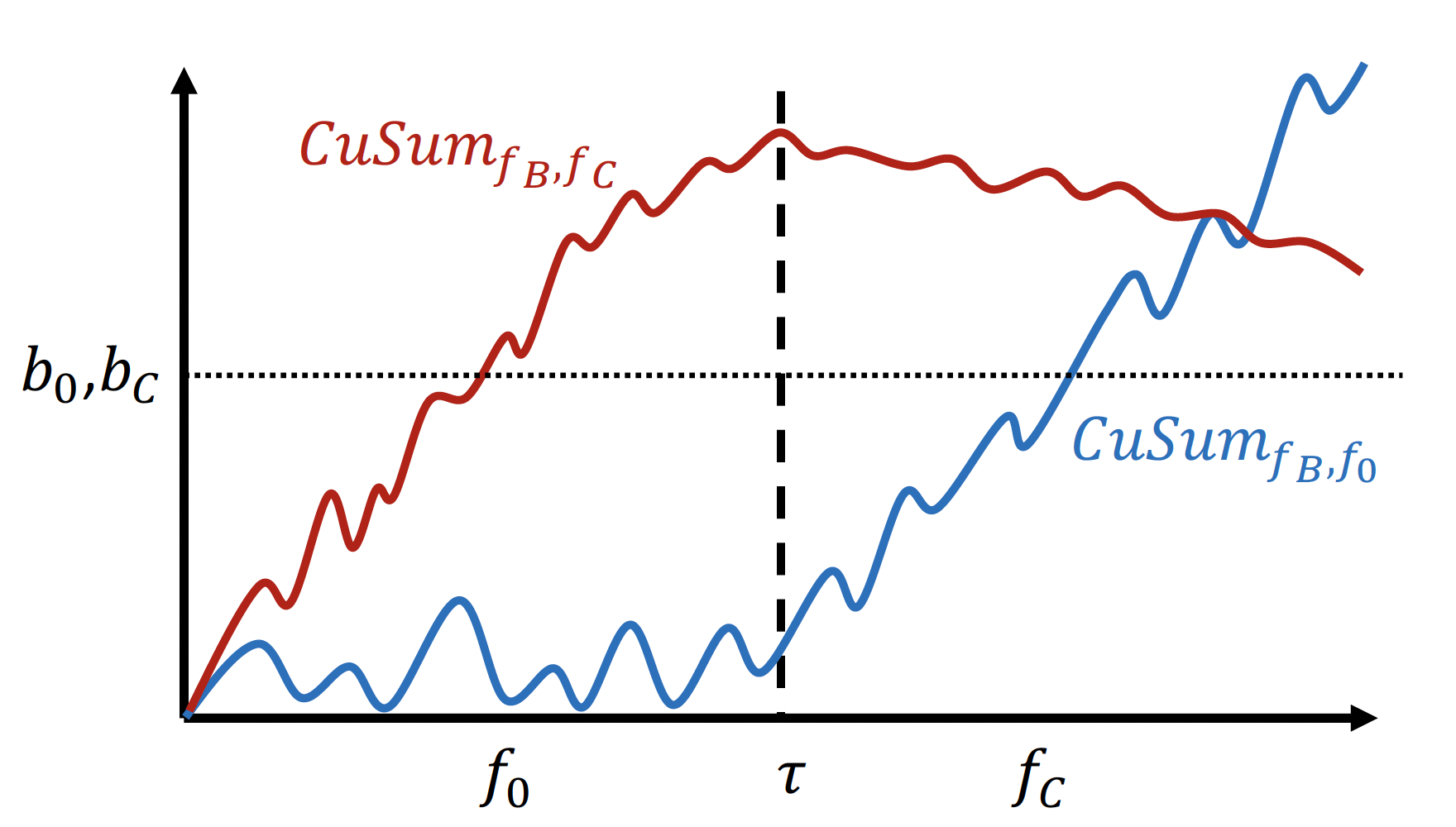}
        \subcaption{Confusing Change}
        \label{fig:nontrivial-confusing-separate}
    \end{minipage}
    \hfill
    \begin{minipage}[b]{0.48\linewidth}
        \includegraphics[width=\columnwidth]{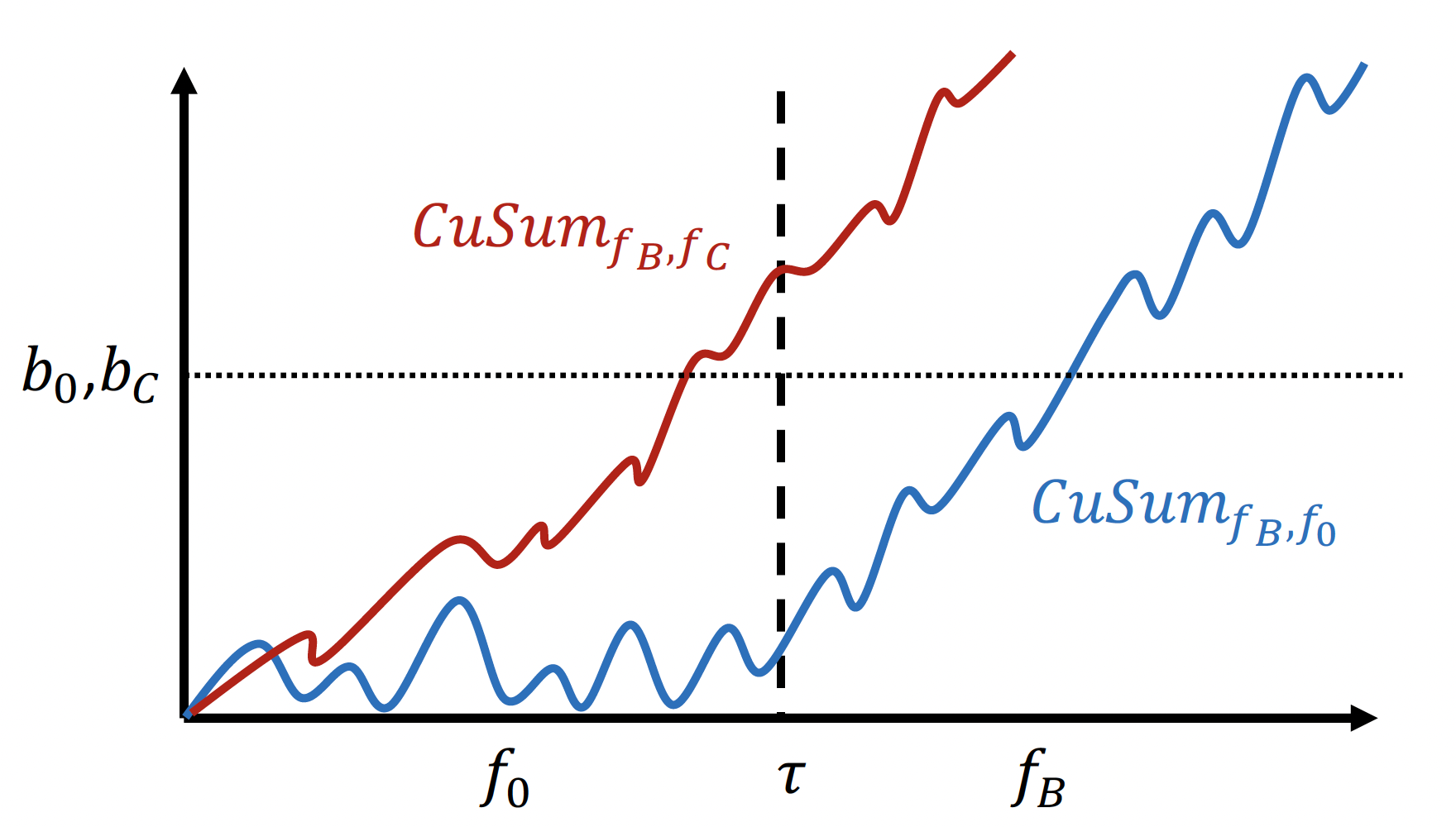}
        \subcaption{Bad Change}
        \label{fig:nontrivial-bad-separate}
    \end{minipage}
    \hfill
    \caption{Illustration of $\text{CuSum}_{f_B, f_0}$ (Blue Lines) and $\text{CuSum}_{f_B, f_C}$ (Red Lines) Behaviors in \textbf{Scenario 3} - $\text{CuSum}_{f_B, f_0}$ fails in distinguishing a bad change from a confusing change, and $\text{CuSum}_{f_B, f_C}$ fails in avoiding raising false alarm during pre-change stage.}
    \label{fig:scenario-3}
\end{figure}

Scenario 3 poses challenges beyond the capabilities of standard single CuSum procedures. Specifically, as illustrated in Figure~\ref{fig:nontrivial-confusing-separate}, in Scenario 3 $\text{CuSum}_{f_B, f_0}$ generally increases under a confusing change distribution and hence is likely to trigger a false alarm. Figure~\ref{fig:scenario-3} also shows that in Scenario 3, the standard $\text{CuSum}_{f_B, f_C}$ procedure is likely to raise a false alarm in a pre-change stage as it has positive drift under $f_0$. Hence, standard single CuSum procedures $\text{CuSum}_{f_B, f_0}$ and $\text{CuSum}_{f_B, f_C}$ cannot solely address this problem. It is worth noting that neither can a naive combination of the standard single CuSum procedures, such as separately launching $\text{CuSum}_{f_B, f_0}$ and $\text{CuSum}_{f_B, f_C}$ simultaneously from the beginning and stopping when both statistics pass thresholds, successfully tackle this problem. Indeed, Figure~\ref{fig:nontrivial-confusing-separate} shows that if the pre-change stage is long, i.e, change point $\nu$ is large, $\text{CuSum}_{f_B, f_C}$ may be much larger than the threshold at the time a confusing change occurs such that it would not fall below the threshold before $\text{CuSum}_{f_B, f_0}$ passes the threshold, and therefore a false alarm would be triggered for a confusing change.


\begin{algorithm}[t]
\small
    \caption{Successive CuSum (\texttt{S-CuSum}) Procedure}\label{alg:S-CuSum}
    \begin{algorithmic}
        \State \textbf{Input:} $f_0, f_C, f_B$, $b_0, b_C$
        \State \textbf{Initialize:} $t = 0$, $\text{CuSum}_{W}[t]\gets 0$, $\text{CuSum}^{\texttt{S}}_{\Lambda}[t] \gets 0$
        \While{1}
            \State $t \gets t + 1$
            \If{$\text{CuSum}_{W}[t-1] < b_0$}
                \State $\text{CuSum}_{W}[t] \gets \left(\text{CuSum}_{W}[t-1] + W[t]\right)^{+}$
            \Else
                \State $\text{CuSum}^{\texttt{S}}_{\Lambda}[t] \gets \left(\text{CuSum}^{\texttt{S}}_{\Lambda}[t-1] + \Lambda[t]\right)^{+}$
                \If{$\text{CuSum}^{\texttt{S}}_{\Lambda}[t-1] \geq b_C$}
                    \State $T_{\texttt{S-CuSum}} \gets t$
                \State Break 
                \EndIf
            \EndIf
        \EndWhile
        \State \textbf{Output:} stopping time $T_{\texttt{S-CuSum}}$
    \end{algorithmic}
\end{algorithm}

\begin{figure}[t]
    \centering
    \begin{minipage}[b]{0.49\linewidth}
        \includegraphics[width=\columnwidth]{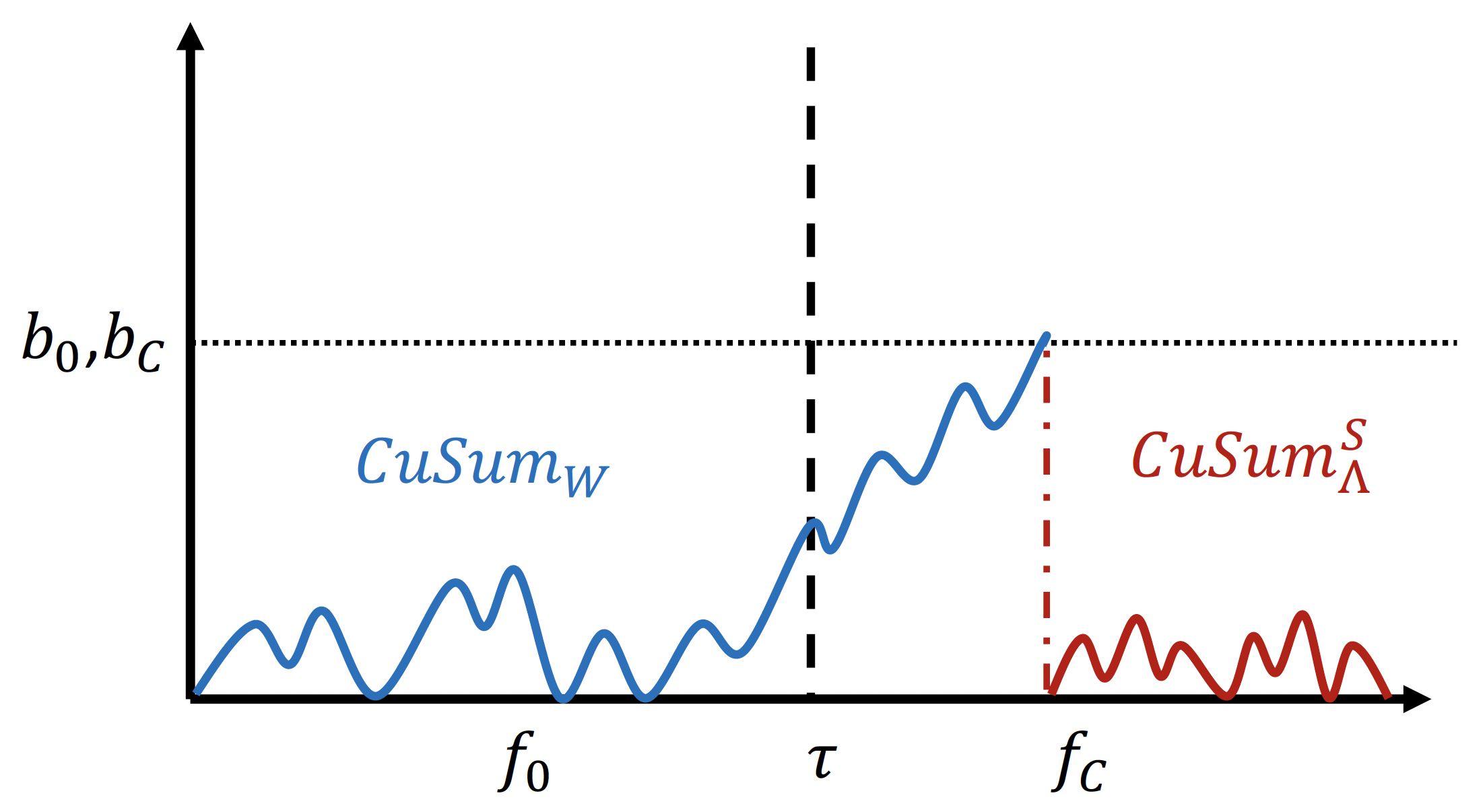}
        \subcaption{Confusing Change}
        \label{fig:nontrivial-confusing-straightforward}
    \end{minipage}
    \hfill
    \begin{minipage}[b]{0.49\linewidth}
        \includegraphics[width=\columnwidth]{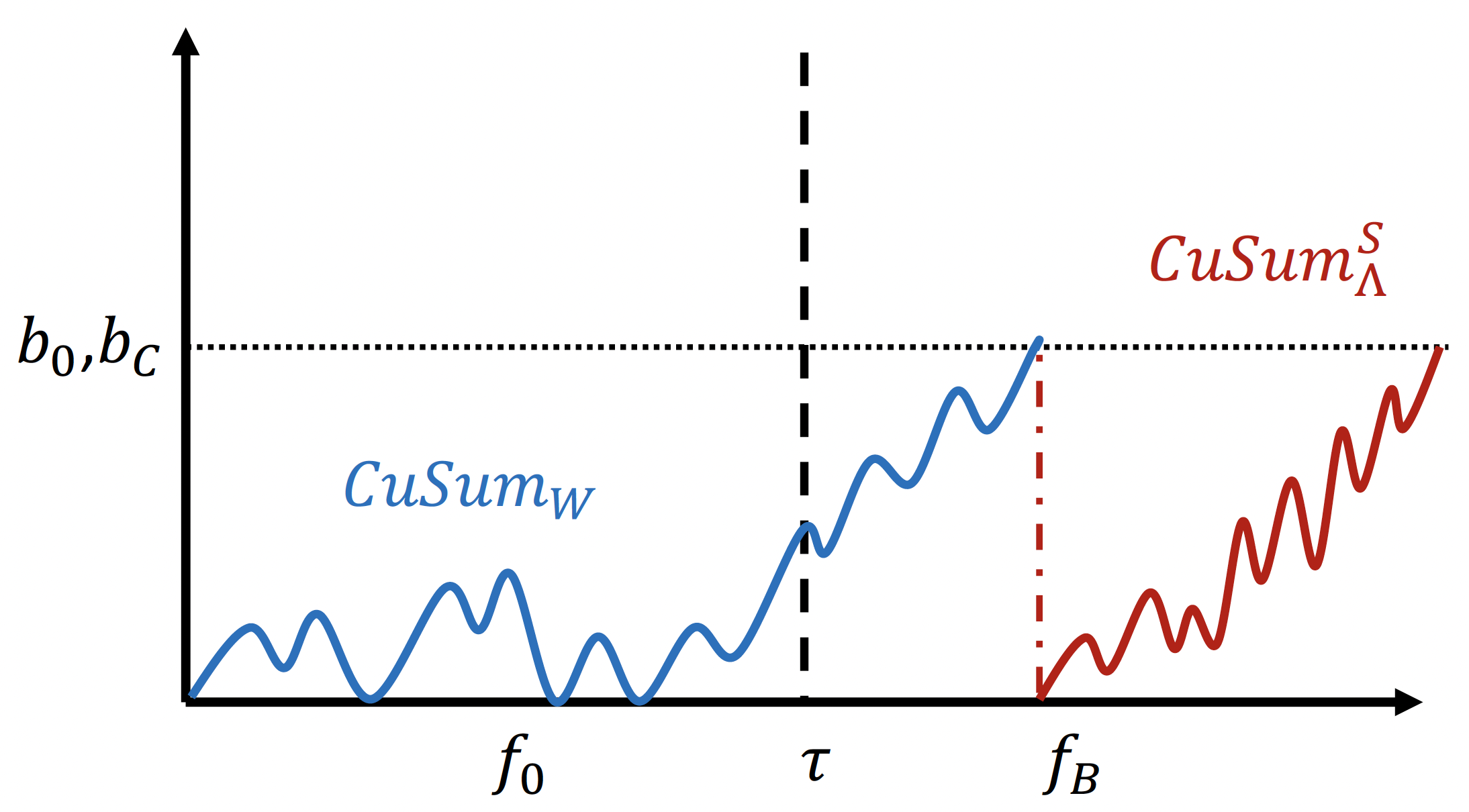}
        \subcaption{Bad Change}
        \label{fig:nontrivial-bad-straightforward}
    \end{minipage}
    \caption{Illustration of \texttt{S-CuSum} Behaviors in Scenario 2 \& 3 - only in the bad change distribution would both $\text{CuSum}_W$ and $\text{CuSum}^{\texttt{S}}_\Lambda$ pass the threshold. Note that \texttt{S-CuSum} also works in Scenario 1.}
    \label{fig:S-CuSum}
\end{figure}

\section{Successive CuSum and Joint CuSum Procedures}\label{sec:algos}

In this section, we propose two new procedures, Successive CuSum (\texttt{S-CuSum}) and Joint CuSum (\texttt{J-CuSum}), that work for all scenarios of the QCD with confusing change problems. We begin by introducing the two hypothesis tests corresponding to this problem. We have one hypothesis test that aims to determine whether we are in the pre-change state or in a post-change state at time $t$:
\begin{align}
    &\mathcal{H}_0[t]: \tau > t\\
    &\mathcal{H}_1[t]: \tau \leq t,
\end{align}
and another hypothesis test that aims to distinguish between the post-change states (confusing change or bad change):
\begin{align}
    &\mathcal{H}_0[t]: X_t \sim f_C(X)\\
    &\mathcal{H}_1[t]: X_t \sim f_B(X).
\end{align}
Intuitively, the hypothesis test to determine pre-/post-change suggests testing against log-likelihood 
\begin{align}
    W[t] := \log \left(\frac{f_B(X_t)}{f_0(X_t)}\right),
\end{align}
and the hypothesis test for distinguishing confusing/bad change suggests testing against log-likelihood 
\begin{align}
    \Lambda[t] := \log \left(\frac{f_B(X_t)}{f_C(X_t)}\right).
\end{align}
Since we are only interested in detecting the \textit{bad change}, we should only raise an alarm when both tests favor the alternative hypotheses. This suggests that synthesizing the two tests is the key to our problem.

\begin{algorithm}[t]
\small
    \caption{Joint CuSum (\texttt{J-CuSum}) Policy}\label{alg:J-CuSum}
    \begin{algorithmic}
        \State \textbf{Input:} $f_0, f_C, f_B$, $b_0, b_C$
        \State \textbf{Initialize:} $t = 0$, $\text{CuSum}_W[t]\gets 0$, $\text{CuSum}^{\texttt{J}}_\Lambda[t] \gets 0$
        \While{1}
            \State $t \gets t + 1$
            \If{$\text{CuSum}_W[t-1] < b_0$}
                \State $\text{CuSum}_W[t] \gets \left(\text{CuSum}_W[t-1] + W[t]\right)^{+}$
            \EndIf
            \If{$\text{CuSum}^{\texttt{J}}_\Lambda[t-1] < b_C$}
                \State $\text{CuSum}^{\texttt{J}}_\Lambda[t] \gets \left(\text{CuSum}^{\texttt{J}}_\Lambda[t-1] + \Lambda[t]\right)^{+}$
            \EndIf
            \If{$\text{CuSum}_W[t] \geq b_0$, $\text{CuSum}^{\texttt{J}}_\Lambda[t] \geq b_C$}
                \State $T_{\texttt{J-CuSum}} \gets t$
                \State Break 
            \ElsIf{$\text{CuSum}_W[t]\leq 0$}
                \State $\text{CuSum}_W[t] \gets 0$
                \State $\text{CuSum}^{\texttt{J}}_\Lambda[t] \gets 0$
            \EndIf
        \EndWhile
        \State \textbf{Output:} stopping time $T_{\texttt{J-CuSum}}$
    \end{algorithmic}
\end{algorithm}

\begin{figure}[t]
    \begin{minipage}[b]{0.48\linewidth}
        \includegraphics[width=\columnwidth]{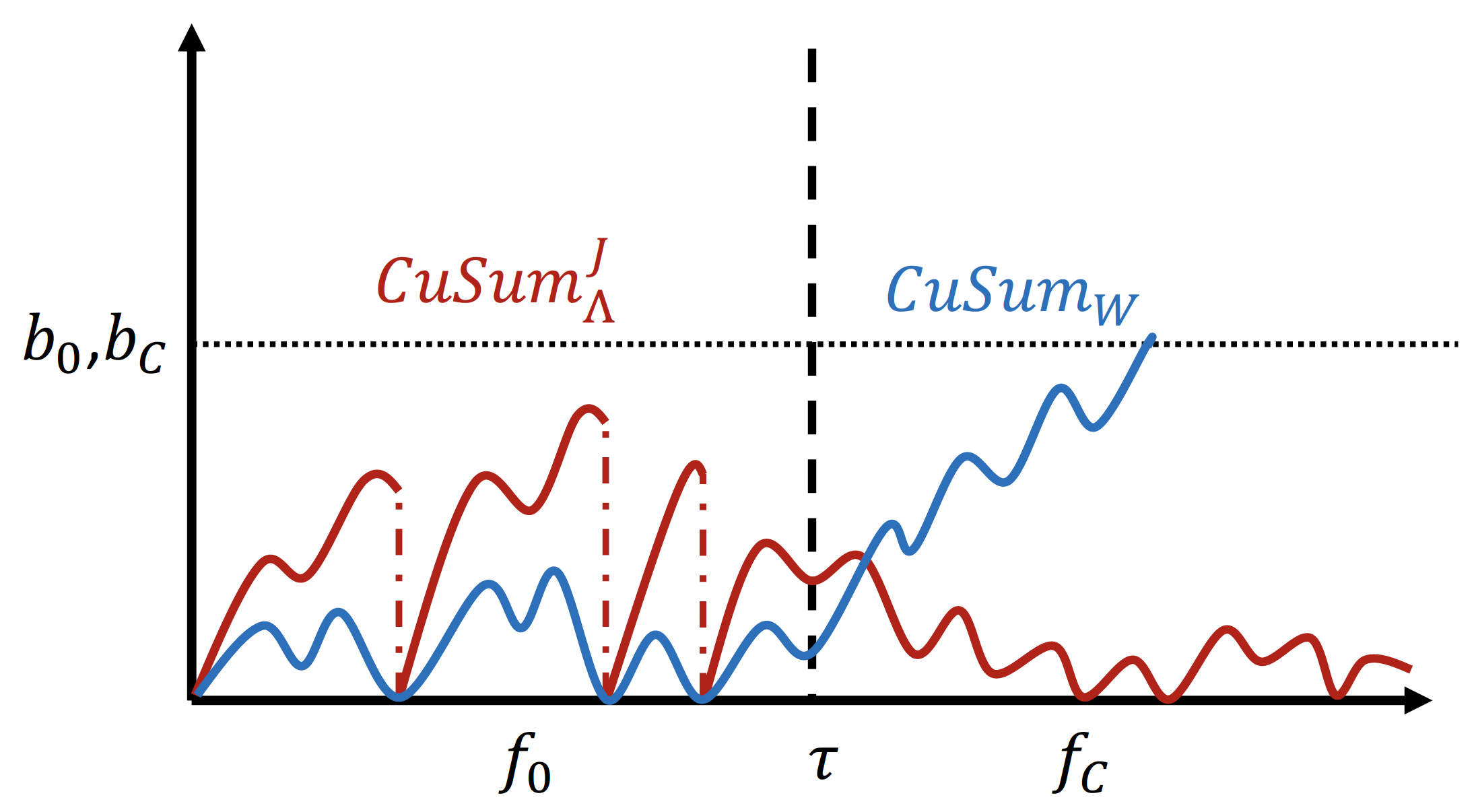}
        \subcaption{Confusing Change}
    \end{minipage}
    \hfill
    \begin{minipage}[b]{0.48\linewidth}
        \includegraphics[width=\columnwidth]{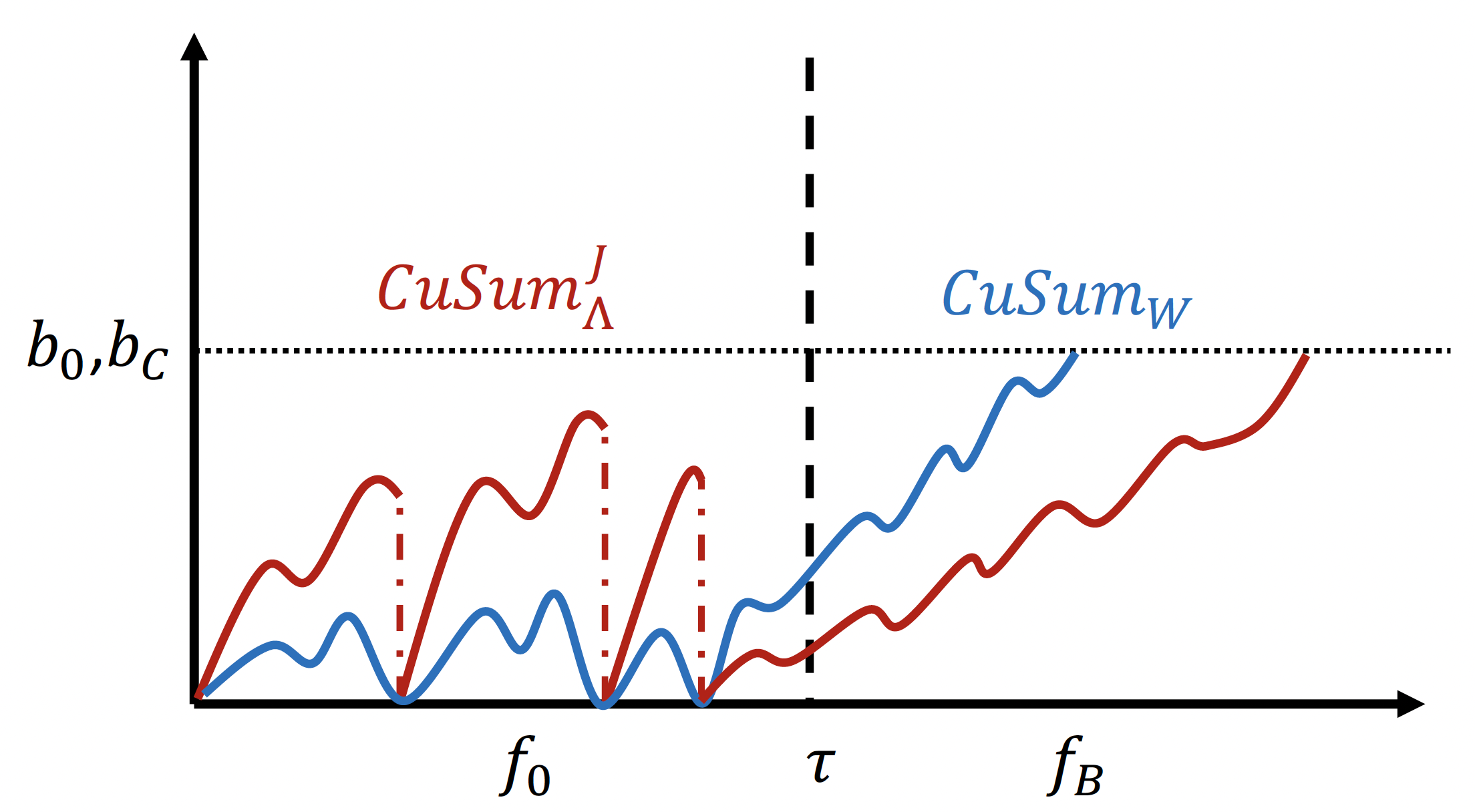}
        \subcaption{Bad Change}
    \end{minipage}
    \caption{Illustration of \texttt{J-CuSum} Behaviors in Scenario 3 - the resetting of $\text{CuSum}^{\texttt{J}}_\Lambda$ (Red Lines) prevents \texttt{J-CuSum} from raising false alarm for confusing change, and applying $\text{CuSum}_W$ (Blue Lines) and $\text{CuSum}^{\texttt{J}}_\Lambda$ simultaneously shortens detection delay. Note that \texttt{J-CuSum} also applies to Scenario 1 \& 2.}
    \label{fig:nontrivial-J-CuSum}
\end{figure}

The core idea of the Successive CuSum (\texttt{S-CuSum}) procedure is to prevent test statistic w.r.t. $\Lambda[t]$ from passing the threshold in the pre-change stage by only launching it after the test statistic w.r.t. $W[t]$ has passed the threshold. Specifically, we let the test statistic w.r.t. $W[t]$ be
\begin{align}\label{eq:CuSumW}
    &\text{CuSum}_W[t]:=\notag\\
    & \begin{cases}
        0, & \text{if } t=0,\\
        \left(\text{CuSum}_W[t-1] + W[t]\right)^{+}, & \text{if $\text{CuSum}_W[t-1] < b_0$, }\\
        \text{CuSum}_W[t-1], & \text{otherwise. }
    \end{cases}
\end{align}
And let the test statistic w.r.t. $\Lambda[t]$ be
\begin{align}
    &\text{CuSum}^{\texttt{S}}_\Lambda[t]:=\notag\\ 
    & \begin{cases}
        0, &\text{if } t=0 \text{ or } \text{CuSum}_W[t] < b_0,\\
        \left(\text{CuSum}^{\texttt{S}}_\Lambda[t-1] + \Lambda[t]\right)^{+}, & \text{otherwise.}
    \end{cases}
\end{align}
That is, \texttt{S-CuSum} stops updating $\text{CuSum}_W$ once it passes the threshold and starts updating $\text{CuSum}^{\texttt{S}}_\Lambda$. An alarm is triggered when $\text{CuSum}^{\texttt{S}}_\Lambda$ passes the threshold, i.e.,
\begin{align}
    T_{\texttt{S-CuSum}} &:= \inf\{t\geq 1: \text{CuSum}^{\texttt{S}}_\Lambda[t] \geq b_C\}\\
    &\equiv\inf\{t\geq 1: \text{CuSum}_W[t] \geq b_0, \text{CuSum}^{\texttt{S}}_\Lambda[t] \geq b_C\}.
\end{align}
Pseudocode for \texttt{S-CuSum} is given in Algorithm~\ref{alg:S-CuSum}. As illustrated in Figure~\ref{fig:S-CuSum}, $\text{CuSum}_W$ most likely passes the threshold after a change has occurred. If the change is a confusing change, $\text{CuSum}^{\texttt{S}}_\Lambda$ always has a negative drift and will rarely pass the threshold; if the change is a bad change, $\text{CuSum}^{\texttt{S}}_\Lambda$ generally increases and passes the threshold quickly.

While \texttt{S-CuSum} effectively detects the bad change and ignores the confusing change, its detection delay leaves room for improvement. Toward this, we propose Joint CuSum (\texttt{J-CuSum}), which incorporates two tests w.r.t. $W[t]$ and $\Lambda[t]$ respectively in a more involved way. Specifically, \texttt{J-CuSum} utilizes $\text{CuSum}_W$ (as defined in \eqref{eq:CuSumW}) and let
\begin{align}
    &\text{CuSum}^{\texttt{J}}_\Lambda[t]:=\notag\\ 
    & \begin{cases}
        0, &\text{ if $t=0$ or $\text{CuSum}_W[t] \leq 0$,}\\
        \left(\text{CuSum}^{\texttt{J}}_\Lambda[t-1] + \Lambda[t]\right)^{+}, & \text{ if $\text{CuSum}^{\texttt{J}}_\Lambda[t-1] < b_C$,}\\
        \text{CuSum}^{\texttt{J}}_\Lambda[t-1],  & \text{otherwise.}
    \end{cases}
\end{align}
That is, \texttt{J-CuSum} prevents $\text{CuSum}^{\texttt{J}}_\Lambda$ from passing the threshold in the pre-change stage by resetting it to zero whenever $\text{CuSum}_W$ hits zero. And \texttt{J-CuSum} also raises an alarm when both statistics pass the threshold, i.e.,
\begin{align}
    T_{\texttt{J-CuSum}} :=\inf \left\{t\geq 1: \text{CuSum}_W[t] \geq b_0, \text{CuSum}^{\texttt{J}}_\Lambda[t] \geq b_C\right\}.
\end{align}
The pseudocode of \texttt{J-CuSum} is presented in Algorithm~\ref{alg:J-CuSum}, and Figure~\ref{fig:nontrivial-J-CuSum} illustrates how \texttt{J-CuSum} works.


\section{Theoretical Guarantees}\label{sec:analysis}

In this section, we discuss the theoretical properties of the QCD with confusing change problems and our proposed procedures. 

We first study the universal lower bound on the $\text{WADD}_{f_B}$ for any procedure whose run length to false alarm is no smaller than $\gamma$.

\begin{theorem}[Universal Detection Delay Lower Bound] \label{thm:universal-lower}
As $\gamma \rightarrow \infty$, we have
\begin{align}
    &\inf_{T\in \mathcal{C}_{\gamma_0, \gamma_C}}\text{WADD}_{f_B}(T) \notag\\
    &\geq \frac{\log(\gamma)(1-o(1))}{\min\{D_{KL}(f_B||f_0), D_{KL}(f_B||f_C)\}+o(1)}.
\end{align}
\end{theorem}

The full proof of Theorem~\ref{thm:universal-lower} is deferred to Appendix~\ref{sec:proof-universal-lower}. The key step in this proof is to utilize our false alarm requirement, \eqref{eq:stopping-time-set}, and the proof-by-contradiction argument as in the proof of~\cite[Theorem 1]{lai1998information}. Intuitively, the detection delay of a stopping rule $T \in \mathcal{C}_\gamma$ not only depends on $D_{KL}(f_B||f_0)$ but also depends on $D_{KL}(f_B||f_C)$ as it is required to have at least $\gamma$ run length to false alarm for confusing change.

In the following, we analyze the run length to false alarm and detection delay of \texttt{S-CuSum}.

\begin{theorem}[\texttt{S-CuSum} False Alarm Lower Bound]\label{thm:SCuSum-lower}
With $b_0 = b_C = \log \gamma$, we have
    \begin{align}
        T_{\texttt{S-CuSum}} \in \mathcal{C}_{\gamma},
    \end{align}
    where $\mathcal{C}_{\gamma}$ is defined in Eq.~\eqref{eq:stopping-time-set}.
\end{theorem}

The full proof of Theorem~\ref{thm:SCuSum-lower} is given in Appendix~\ref{sec:proof-SCuSum-lower}. Note the \texttt{S-CuSum} only triggers an alarm when both $\text{CuSum}_W$ and $\text{CuSum}^{\texttt{S}}_\Lambda$ pass their thresholds. Hence, to show that $T_{\texttt{S-CuSum}} \in \mathcal{C}_{\gamma}$, we need to show that $\mathbb{E}_{\infty}[T_{\text{CuSum}_W}] \geq \gamma$ and that $\mathbb{E}_{1, f_C}[T_{\text{CuSum}_{f_B, f_C}}] \geq \gamma$. We establish both inequalities by relating the CuSum statistics to their corresponding Shiryaev-Roberts statistics and utilizing the martingale properties of the Shiryaev-Roberts statistics~\cite{veeravalli2014quickest}.

\begin{theorem}[\texttt{S-CuSum} Detection Delay Upper Bound]\label{thm:SCuSum-upper}
With $b_0 = b_C = \log \gamma$, as $\gamma \rightarrow \infty$, we have
\begin{align}
    &\text{WADD}_{f_B}(T_{\texttt{S-CuSum}})\notag\\
    &\leq \left(\frac{\log \gamma}{D_{KL}(f_B||f_0)}+\frac{\log \gamma}{D_{KL}(f_B||f_C)}\right)(1+o(1))\\
    &\leq \frac{2\log\gamma}{\min\{D_{KL}(f_B||f_0), D_{KL}(f_B||f_C)\}}(1+o(1)).
\end{align}
\end{theorem}

The full proof of Theorem~\ref{thm:SCuSum-upper} is given in Appendix~\ref{sec:proof-SCuSum-upper}. Because $\text{CuSum}_W[t]$ and $\text{CuSum}^{\texttt{S}}_\Lambda[t]$ are always non-negative and are zero when $t=0$, the worse-case average detection delay occurs when the change point $\nu = 1$. And by the algorithmic property of \texttt{S-CuSum}, $\mathbb{E}_{1, f_B}[T_{\texttt{S-CuSum}}]$ equals the sum of $\mathbb{E}_{1, f_B}[T_{\text{CuSum}_W}]$ and $\mathbb{E}_{1, f_B}[T_{\text{CuSum}_{f_B, f_C}}]$. We upper bound both using a generalized Weak Law of Large Numbers, Lemma~\ref{lemma:WLLN}.

\begin{figure*}[t]
    \begin{minipage}[b]{0.3\linewidth}
        \includegraphics[width=\columnwidth]{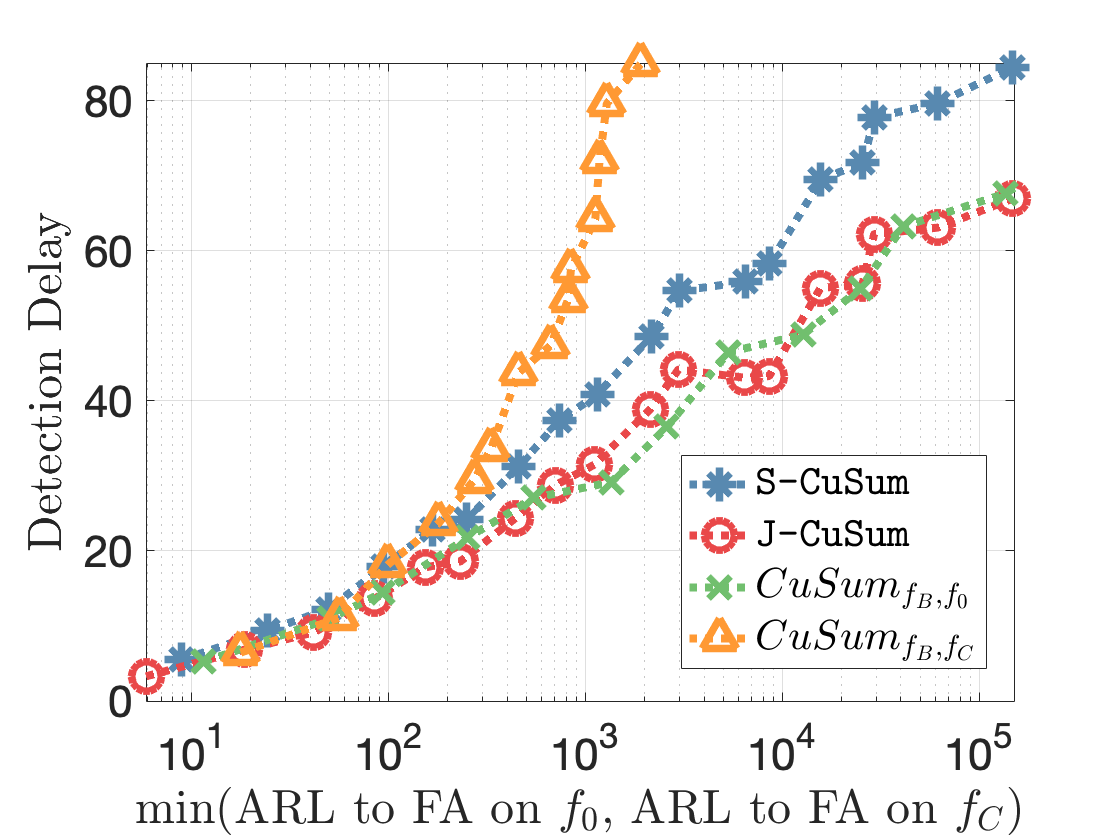} 
        \subcaption{\textbf{Scenario 1}}
        \label{fig:S1}
    \end{minipage}
    \hfill
    \begin{minipage}[b]{0.3\linewidth}
        \includegraphics[width=\columnwidth]{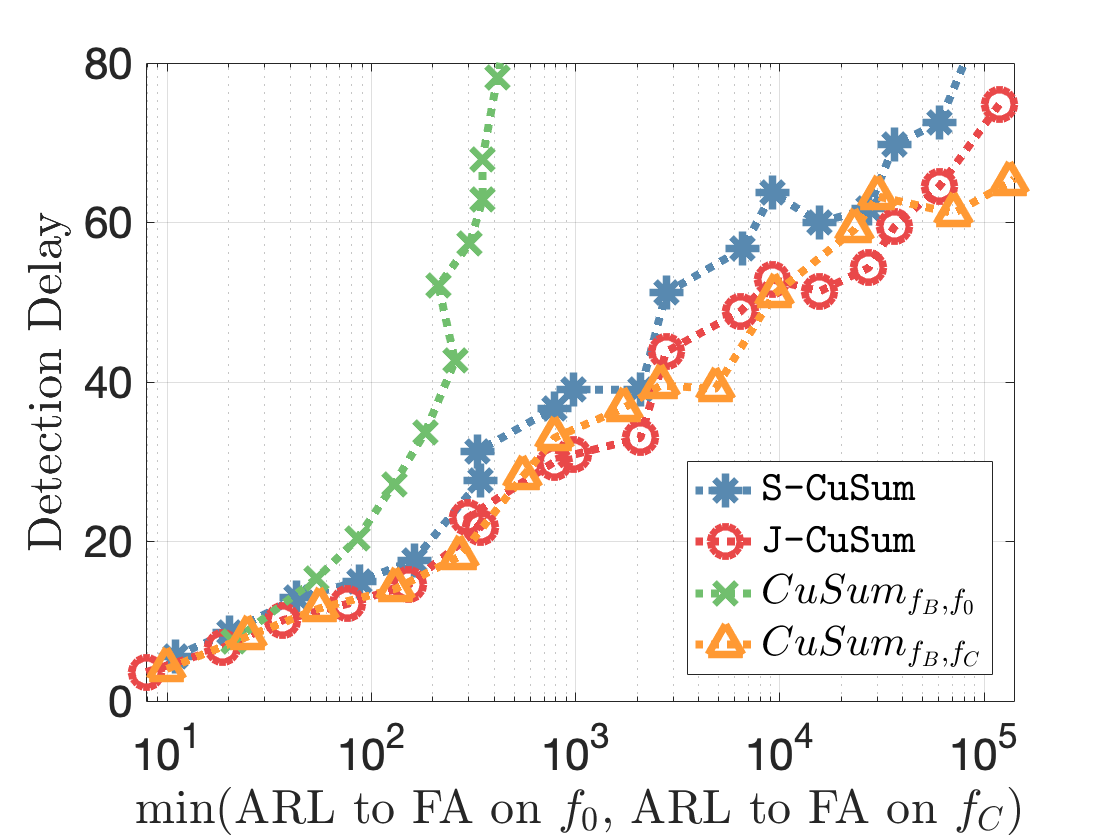}
        \subcaption{\textbf{Scenario 2}}
        \label{fig:S2}
    \end{minipage}
    \hfill
    \begin{minipage}[b]{0.3\linewidth}
        \includegraphics[width=\columnwidth]{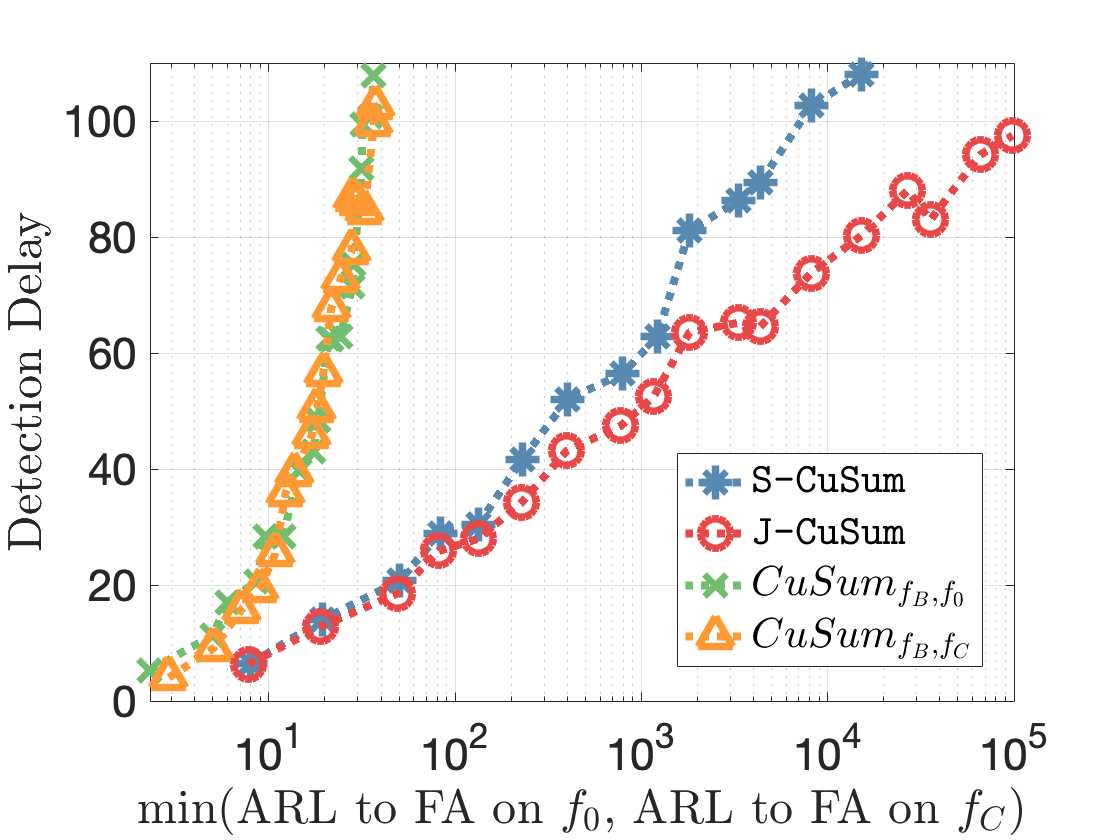}
        \subcaption{\textbf{Scenario 3}}
        \label{fig:S3}
    \end{minipage}
    \caption{Numerical Comparison between \texttt{S-CuSum} and \texttt{J-CuSum} and baselines, $\text{CuSum}_{f_B, f_0}$ and $\text{CuSum}_{f_B, f_C}$ in the three scenarios specified in Section~\ref{sec:scenario}.}
    \label{fig:simulations}
\end{figure*}

In the following, we analyze the run length to false alarm and detection delay of \texttt{J-CuSum}.

\begin{theorem}[\texttt{J-CuSum} False Alarm Lower Bound]\label{thm:JCuSum-lower}
With $b_0 = b_C = \log \gamma$, we have
    \begin{align}
        T_{\texttt{J-CuSum}} \in \mathcal{C}_{\gamma},
    \end{align}
    where $\mathcal{C}_{\gamma}$ is defined in Eq.~\eqref{eq:stopping-time-set}.
\end{theorem}

The full proof of Theorem~\ref{thm:JCuSum-lower} is deferred to Appendix~\ref{sec:proof-JCuSum-lower}. By the algorithmic property of $\texttt{J-CuSum}$, when confusing change occurs and the change point is just right before $\text{CuSum}_{\Lambda}^{\texttt{J}}$ being reset by $\text{CuSum}_W$, $\texttt{J-CuSum}$ has the shortest average run time for $\text{CuSum}_\Lambda^{\texttt{J}}$ to pass the threshold. To analyze the run length between a reset to the next reset, we follow~\cite{siegmund1985sequential} to define the stopping time of $\text{CuSum}_W$ in terms of the stopping times of a sequence of sequential probability ratio tests. It follows from the approximation of the stopping time of a corresponding sequential probability ratio test given in~\cite{siegmund1985sequential} that the expected value of $\text{CuSum}_\Lambda^{\texttt{J}}$ just right before a reset is almost zero. Therefore, we can utilize results in the proof of Theorem~\ref{thm:SCuSum-lower} to prove that $T_{\texttt{J-CuSum}} \in \mathcal{C}_{\gamma}$.

\begin{theorem}[\texttt{J-CuSum} Detection Delay Upper Bound]\label{thm:JCuSum-upper}
With $b_0 = b_C = \log \gamma$, as $\gamma \rightarrow \infty$, we have
\begin{align}
    &\text{WADD}_{f_B}(T_{\texttt{J-CuSum}}) \notag\\
    &\lessapprox \left(\frac{\log \gamma}{D_{KL}(f_B||f_0)}+\frac{\log \gamma}{D_{KL}(f_B||f_C)}\right)(1+o(1))\\
    &\leq \frac{2\log\gamma}{\min\{D_{KL}(f_B||f_0), D_{KL}(f_B||f_C)\}}(1+o(1)).
\end{align}
\end{theorem}

The full proof of Theorem~\ref{thm:JCuSum-upper} is deferred to Appendix~\ref{sec:proof-JCuSum-upper}. 
Because $\text{CuSum}_W[t]$ and $\text{CuSum}^{\texttt{J}}_\Lambda[t]$ are always non-negative and are zero when $t=0$, the worst-case average detection delay occurs when the change point $\nu = 1$. 
And by the algorithmic property of \texttt{J-CuSum}, $\mathbb{E}_{1, f_B}[T_{\texttt{J-CuSum}}]$ is upper bounded by both $\mathbb{E}_{1, f_B}[T_{\text{CuSum}_W}]$ and $\mathbb{E}_{1, f_B}[T_{\text{CuSum}^{\texttt{J}}_\Lambda}]$. We upper bound $\mathbb{E}_{1, f_B}[T_{\text{CuSum}_W}]$ using Lemma~\ref{lemma:WLLN}. As for $\mathbb{E}_{1, f_B}[T_{\text{CuSum}^{\texttt{J}}_\Lambda}]$, we follow~\cite{siegmund1985sequential} to define the stopping time of $\text{CuSum}_W$ in terms of the stopping times of a sequence of sequential probability ratio tests and approximate the run length to the last resetting.


\section{Numerical Results}\label{sec:simulations}

In this section, we numerically compare \texttt{S-CuSum} and \texttt{J-CuSum} to baselines, $\text{CuSum}_{f_B, f_0}$ and $\text{CuSum}_{f_B, f_C}$. Specifically, we conduct simulations in each of the three scenarios discussed in Section~\ref{sec:scenario}:
\begin{itemize}
    \item for \textbf{Scenario 1}, we let $f_0 = \mathcal{N}(0, 1)$, $f_C = \mathcal{N}(-0.5, 1)$, and $f_B = \mathcal{N}(0.5, 1)$;
    \item for \textbf{Scenario 2}, we let $f_0 = \mathcal{N}(0, 1)$, $f_C = \mathcal{N}(0.7, 1)$, and $f_B = \mathcal{N}(1.2, 1)$;
    \item for \textbf{Scenario 3}, we let $f_0 = \mathcal{N}(0, 1)$, $f_C = \mathcal{N}(1, 1)$, and $f_B = \mathcal{N}(0.5, 1)$.
\end{itemize}
For each scenario, we run all procedures under $\mathbb{P}_{1, f_B}$, $\mathbb{P}_{\infty}$, and $\mathbb{P}_{1, f_C}$ with varying thresholds (for \texttt{S-CuSum} and \texttt{J-CuSum}, we let $b_0=b_C$) to learn these procedures' detection delays, pre-change run lengths to false alarms, and run lengths to false alarm for confusing change respectively. For each scenario, underlying distribution, and threshold, we perform $60$ independent trials and report the average. 

In Figure~\ref{fig:simulations}, we plot the average detection delays of the procedures against their average run lengths to false alarm for pre-change or for confusing change, whichever is smaller.
This is because our false alarm requirement, Eq.~\eqref{eq:stopping-time-set}, asks for both the run lengths to false alarm for pre-change and confusing change to be no smaller than the same threshold; hence, we plot against whichever is smaller to make sure that the requirement is fulfilled. Moreover, We plot run lengths to false alarm on a logarithmic scale while plotting detection delay on a linear scale. Hence, straight lines on the figures indicate that detection delays grow logarithmically with regard to run lengths to false alarm; whereas the steeply rising lines on the figures indicate that detection delays grow super-logarithmically with regard to run lengths to false alarm.

First, Figures~\ref{fig:S1}, \ref{fig:S2}, and \ref{fig:S3} respectively corroborate our discussion in Section~\ref{sec:scenario} that $\text{CuSum}_{f_B, f_0}$ suffices in \textbf{Scenario 1}, $\text{CuSum}_{f_B, f_C}$ suffices in \textbf{Scenario 2}, but neither a single $\text{CuSum}_{f_B, f_0}$ nor a single $\text{CuSum}_{f_B, f_C}$ suffice in \textbf{Scenario 3} to detect a bad change quickly while avoiding raising false alarm for a confusing change. Indeed, Figure~\ref{fig:S3} shows that both $\text{CuSum}_{f_B, f_0}$ and $\text{CuSum}_{f_B, f_C}$ incur very large average detection delays in order to achieve same average run lengths to false alarm as \texttt{S-CuSum} or \texttt{J-CuSUm}.

Second, Figures~\ref{fig:S1}, \ref{fig:S2}, and \ref{fig:S3} show that both \texttt{S-CuSum} and \texttt{J-CuSum} perform well in all three possible scenarios, namely under all kinds of pre-change, bad change, and confusing change distributions. Indeed, in Figure~\ref{fig:simulations}, the lines of \texttt{S-CuSum} and \texttt{J-CuSum} are straight in all three graphs, meaning that their detection delays grow logarithmically with regard to their run lengths to false alarm.

Finally, as shown in Figures~\ref{fig:S1}, \ref{fig:S2}, and \ref{fig:S3}, our simulations empirically support our discussion in Section~\ref{sec:algos} that the detection delay of \texttt{S-CuSum} still has room to improve and \texttt{J-CuSum} reduces the detection delay by allowing the $\text{CuSum}^{\texttt{J}}_\Lambda$ to be launched earlier than $\text{CuSum}^{\texttt{S}}_\Lambda$ but still avoiding raising false alarms.

\section{Conclusion and Discussion}\label{sec:conclusion}

In this paper, we investigated a quickest change detection problem where an initially in-control system can transition into an out-of-control state due to either a bad event or a confusing event. Our goal was to detect the change as soon as possible if a bad event occurs while avoiding raising an alarm if a confusing event occurs. We found that when both 1) the KL-divergence between confusing change distribution $f_C$ and pre-change distribution $f_0$ is larger than that between $f_C$ and bad change distribution $f_B$, and 2) the KL-divergence between $f_0$ and $f_C$ is greater than that between $f_0$ and $f_B$ occurs, typical procedures based on $\text{CuSum}_{f_B, f_0}$ and $\text{CuSum}_{f_B, f_C}$ fail to achieve our objective. Hence, we proposed two new detection procedures \texttt{S-CuSum} and \texttt{J-CuSum} that achieve our objective in all scenarios and provide theoretical guarantees as well as numerical corroborations.  

While the detection delay upper bound of \texttt{J-CuSum} that we obtained is the same as that obtained for \texttt{S-CuSum}, intuitively,  \texttt{J-CuSum} should produce smaller detection delays than those \texttt{S-CuSum} would incur. Indeed, in all our simulations, \texttt{J-CuSum} has smaller detection delays than \texttt{S-CuSum}. We leave closing the theoretical gap between the detection delay upper bound of \texttt{J-CuSum} and the universal lower bound of detection delay for future work.

\section*{Acknowledgement}

The authors would like to thank Lance Kaplan for invaluable insights and engaging discussions throughout the development of this paper.

\appendices
\bibliography{ref}
\bibliographystyle{IEEEtran}


\section{proof of Theorem~\ref{thm:universal-lower}}\label{sec:proof-universal-lower}

\begin{proof}
    We first recall the following generalized version of the Weak Law of Large Number. 
    \begin{lemma}[Lemma A.1 in \cite{fellouris2017multichannel}]\label{lemma:WLLN}
        Let $\{Y_t, t\in \mathbb{N}\}$ be a sequence of random variables i.i.d. on $(\Omega, \mathcal{F}, \mathbb{P})$ with $\mathbb{E}[Y_t] = \mu > 0$, then for any $\epsilon > 0$, as $n\rightarrow \infty$, 
        \begin{align}
            \mathbb{P}\left[\frac{\max_{1\leq k \leq n} \sum_{t=1}^k Y_t}{n} - \mu > \epsilon\right] \rightarrow 0.
        \end{align}
    \end{lemma}    

    Note that
    \begin{align}
        \text{WADD}_{f_B}(T) &:= \sup\limits_{\nu \geq 1}\mathbb{E}_{\nu, f_B}[T-\nu|T \geq \nu]\\
        &\geq \mathbb{E}_{\nu, f_B}[T-\nu|T \geq \nu]\\
        &\overset{(a)}{\geq} \mathbb{P}_{\nu, f_B}[T-\nu\geq \alpha_\gamma|T\geq \nu]\times \alpha_\gamma,
    \end{align}
    where inequality (a) is by the Markov inequality.
    It then suffices to show that as $\gamma \rightarrow \infty$,
    \begin{align}
        \mathbb{P}_{\nu, f_B}[T-\nu\geq \alpha_\gamma|T\geq \nu] \rightarrow 1,
    \end{align}
    or equivalently, 
    \begin{align}\label{eq:lower-proof-goal}
        \mathbb{P}_{\nu, f_B}[\nu \leq T < \nu +\alpha_\gamma | \tau \geq \nu] \rightarrow 0.
    \end{align}

    We will first show that Eq.~\eqref{eq:lower-proof-goal} holds when
    \begin{align}\label{eq:lower-proof-part-1-begin}
        &\alpha_\gamma = \frac{1}{D_{KL}(f_B||f_0)+\epsilon}\log \gamma^{(1-\epsilon)}, & \epsilon >0,
    \end{align}
    using a change-of-measure argument.
    Specifically, 
    \begin{align}
        &\mathbb{P}_{f_0}[\nu \leq T < \nu+\alpha_\gamma]\notag\\
        &= \mathbb{E}_{f_0}[\mathds{1}_{\{\nu \leq T < \nu+\alpha_\gamma\}}]\\
        &\overset{(a)}{=} \mathbb{E}_{\nu, f_B}\Big[\mathds{1}_{\{\nu \leq T < \nu+\alpha_\gamma\}} \frac{\mathbb{P}_{f_0}(H_T)}{\mathbb{P}_{\nu, f_B}(H_T)}\Big]\\
        &\geq \mathbb{E}_{\nu, f_B}\Big[\mathds{1}_{\{\nu \leq T < \nu+\alpha_\gamma, \log(\frac{\mathbb{P}_{f_0}(H_T)}{\mathbb{P}_{\nu, f_B}(H_T)}) \geq -a\}} \frac{\mathbb{P}_{f_0}(H_T)}{\mathbb{P}_{\nu, f_B}(H_T)}\Big]\\
        &\geq e^{-a} \mathbb{P}_{\nu, f_B}\Big[\nu \leq T < \nu+\alpha_\gamma, \log\Big(\frac{\mathbb{P}_{f_0}(H_T)}{\mathbb{P}_{\nu, f_B}(H_T)}\Big) \geq -a\Big]\\
        &=e^{-a}\mathbb{P}_{\nu, f_B}\Big[\nu \leq T < \nu+\alpha_\gamma, \log\Big(\frac{\mathbb{P}_{\nu, f_B}(H_T)}{\mathbb{P}_{f_0}(H_T)}\Big) \leq a\Big]\\
        &\geq e^{-a} \mathbb{P}_{\nu, f_B}\Big[\nu \leq T < \nu+\alpha_\gamma,\notag\\
        &\qquad\qquad\qquad\qquad \max\limits_{\nu\leq j< \nu+\alpha_\gamma}\log\Big(\frac{\mathbb{P}_{\nu, f_B}(H_j)}{\mathbb{P}_{f_0}(H_j)}\Big) \leq a\Big]\\
        &\overset{(b)}{\geq} e^{-a} \mathbb{P}_{\nu, f_B}\Big[\nu \leq T < \nu+\alpha_\gamma\Big] \notag\\
        &\quad - e^{-a}\mathbb{P}_{\nu, f_B}\Big[ \max\limits_{\nu\leq j< \nu+\alpha_\gamma}\log\Big(\frac{\mathbb{P}_{\nu, f_B}(H_j)}{\mathbb{P}_{f_0}(H_j)}\Big) > a\Big]\label{eq:lower-proof-break}
    \end{align}
    where $H_t = (X_1, ..., X_{t}), t\in \mathbb{N}$; change-of-measure argument (a) holds because $\mathbb{P}_{f_0}$ and $\mathbb{P}_{\nu, f_B}$ are measures over a common measurable space, $\mathbb{P}_{\nu, f_B}$ is $\sigma$-finite, and $\mathbb{P}_{f_0} \ll \mathbb{P}_{\nu, f_B}$; $a$ will be specified later; and inequality (b) is because, for any event $A$ and $B$, $\mathbb{P}[A \cap B] \geq \mathbb{P}[A] - \mathbb{P}[B^c]$.

    The event $\{T\geq \nu\}$ only depends on $H_{\nu-1}$, which follows the same distribution under $\mathbb{P}_{f_0}$ and $\mathbb{P}_{\nu, f_B}$. This implies
    \begin{align}\label{eq:lower-proof-condition}
        \mathbb{P}_{f_0}[T \geq \nu] = \mathbb{P}_{\nu, f_B}[T \geq \nu].
    \end{align}
    By Eq.~\eqref{eq:lower-proof-condition} and reordering Eq.~\eqref{eq:lower-proof-break}, it follows that
    \begin{align}
        &\mathbb{P}_{\nu, f_B}[\nu\leq T< \nu+\alpha_\gamma | T \geq \nu] \notag\\
        &\leq e^{a} \mathbb{P}_{f_0}[\nu \leq T< \nu+\alpha_\gamma | T\geq \nu]\notag\\
        &\quad+\mathbb{P}_{\nu, f_B}\Big[ \max\limits_{\nu\leq j< \nu+\alpha_\gamma}\log\Big(\frac{\mathbb{P}_{\nu, f_B}(H_j)}{\mathbb{P}_{f_0}(H_j)}\Big) > a\Big|T\geq \nu\Big].\label{eq:lower-proof-break-2}
    \end{align}

    To show that the first term at the right-hand side of Eq.~\eqref{eq:lower-proof-break-2} converges to $0$ as $\gamma \rightarrow \infty$, we can utilize the proof-by-contradiction argument as in the proof of~\cite[Theorem 1]{lai1998information}. Let $\alpha_\gamma$ be a positive integer and $\alpha_\gamma < \gamma$. For any $T \in \mathcal{C}_\gamma$, we have $\mathbb{E}_{f_0}[T]\geq \gamma$ and then for some $\nu \geq 1$, $\mathbb{P}_{f_0}[T \geq \nu] > 0$ and 
    \begin{align}
        \mathbb{P}_{f_0}[T\leq \nu +\alpha_\gamma|T\geq \nu] \leq \frac{\alpha_\gamma}{\gamma},
    \end{align}
    because otherwise $\mathbb{P}_{f_0}[T \geq \nu +\alpha_\gamma | T \geq \nu] < 1-\nicefrac{\alpha_\gamma}{\gamma}$ for all $\nu \geq 1$ with $\mathbb{P}_{f_0}[T\geq \nu] > 0$, implying that $\mathbb{E}_{f_0}[T] \leq \gamma$.

    Let $a = \log\gamma^{(1-\epsilon)}$, then
    \begin{align}
        &e^{a} \mathbb{P}_{f_0}[\nu \leq T< \nu+\alpha_\gamma | T\geq \nu] \notag\\
        &\leq e^{a} \frac{\alpha_\gamma}{\gamma} = \frac{\alpha_\gamma}{\gamma^\epsilon} \rightarrow 0, \text{ as } \gamma \rightarrow \infty.
    \end{align}

    We then show that the second term at the right-hand side of Eq.~\eqref{eq:lower-proof-break-2} also converges to $0$ as $\gamma \rightarrow \infty$.
    \begin{align}
        &\mathbb{P}_{\nu, f_B}\Big[\max\limits_{\nu\leq j< \nu+\alpha_\gamma}\log\Big(\frac{\mathbb{P}_{\nu, f_B}(H_j)}{\mathbb{P}_{f_0}(H_j)}\Big) > a\Big|T\geq \nu\Big]\notag\\
        &= \mathbb{P}_{\nu, f_B}\Big[\max\limits_{\nu\leq j< \nu+\alpha_\gamma}\sum_{i=\nu}^j \log\Big(\frac{f_B(X_i)}{f_0(X_i)}\Big) > a |T\geq \nu\Big]\\
        &\overset{(a)}{=}\mathbb{P}_{\nu, f_B}\Big[\max\limits_{\nu\leq j< \nu+\alpha_\gamma}\sum_{i=\nu}^j \log\Big(\frac{f_B(X_i)}{f_0(X_i)}\Big)> a\Big]\\
        &\overset{(b)}{\leq} \mathbb{P}_{\nu, f_B}\Big[\max\limits_{\nu\leq j< \nu+\alpha_\gamma}\sum_{i=\nu}^j \log\Big(\frac{f_B(X_i)}{f_0(X_i)}\Big)\notag\\
        &\qquad\qquad\qquad\qquad> \alpha_\gamma \ (D_{KL}(f_B||f_0)+\epsilon)\Big] \\
        &\rightarrow 0, \text{ as } \gamma \rightarrow \infty. \label{eq:lower-proof-part-1-end}
    \end{align}
    where equality (a) is due to the fact that the event $\{T \geq \nu\}$ is independent from $X_i, \forall i \geq \nu$; inequality (b) is because 
    $a \geq \alpha_\gamma \times (D_{KL}(f_B||f_0)+\epsilon)$;
    and the last step is by applying Lemma~\ref{lemma:WLLN}.

    Similarly, we can show that Eq.~\eqref{eq:lower-proof-goal} holds when
    \begin{align}\label{eq:lower-proof-part-2-begin}
        &\alpha_\gamma = \frac{1}{D_{KL}(f_B||f_C)+\epsilon}\log \gamma^{(1-\epsilon)}, & \epsilon >0,
    \end{align}
    using a change-of-measure argument.
    Specifically, 
    \begin{align}
        &\mathbb{P}_{\nu, f_C}[\nu \leq T < \nu+\alpha_\gamma]\notag\\
        &\overset{(a)}{=} \mathbb{E}_{\nu, f_B}\Big[\mathds{1}_{\{\nu \leq T < \nu+\alpha_\gamma\}} \frac{\mathbb{P}_{\nu, f_C}(H_T)}{\mathbb{P}_{\nu, f_B}(H_T)}\Big]\\
        &\overset{(b)}{\geq} e^{-a} \mathbb{P}_{\nu, f_B}\Big[\nu \leq T < \nu+\alpha_\gamma\Big] \notag\\
        &\quad - e^{-a}\mathbb{P}_{\nu, f_B}\Big[ \max\limits_{\nu\leq j< \nu+\alpha_\gamma}\log\Big(\frac{\mathbb{P}_{\nu, f_B}(H_j)}{\mathbb{P}_{\nu, f_C}(H_j)}\Big) > a\Big]\label{eq:lower-proof-break-3}
    \end{align}
    where $H_t = (X_1, ..., X_{t}), t\in \mathbb{N}$; change-of-measure argument (a) holds because $\mathbb{P}_{\nu, f_C}$ and $\mathbb{P}_{\nu, f_B}$ are measures over a common measurable space, $\mathbb{P}_{\nu, f_B}$ is $\sigma$-finite, and $\mathbb{P}_{\nu, f_C} \ll \mathbb{P}_{\nu, f_B}$; $a$ will be specified later; and inequality (b) is because, for any event $A$ and $B$, $\mathbb{P}[A \cap B] \geq \mathbb{P}[A] - \mathbb{P}[B^c]$.

    The event $\{T\geq \nu\}$ only depends on $H_{\nu-1}$, which follows the same distribution under $\mathbb{P}_{\nu, f_C}$ and $\mathbb{P}_{\nu, f_B}$. This implies
    \begin{align}\label{eq:lower-proof-condition-2}
        \mathbb{P}_{\nu, f_C}[T \geq \nu] = \mathbb{P}_{\nu, f_B}[T \geq \nu].
    \end{align}
    By Eq.~\eqref{eq:lower-proof-condition-2} and reordering Eq.~\eqref{eq:lower-proof-break-3}, it follows that
    \begin{align}
        &\mathbb{P}_{\nu, f_B}[\nu\leq T< \nu+\alpha_\gamma | T \geq \nu] \notag\\
        &\leq e^{a} \mathbb{P}_{\nu, f_C}[\nu \leq T< \nu+\alpha_\gamma | T\geq \nu]\notag\\
        &\quad+\mathbb{P}_{\nu, f_B}\Big[ \max\limits_{\nu\leq j< \nu+\alpha_\gamma}\log\Big(\frac{\mathbb{P}_{\nu, f_B}(H_j)}{\mathbb{P}_{\nu, f_C}(H_j)}\Big) > a\Big|T\geq \nu\Big].\label{eq:lower-proof-break-4}
    \end{align}

    To show that the first term at the right-hand side of Eq.~\eqref{eq:lower-proof-break-4} converges to $0$ as $\gamma \rightarrow \infty$, we can utilize the proof-by-contradiction argument as in the proof of~\cite[Theorem 1]{lai1998information}. Let $\alpha_\gamma$ be a positive integer and $\alpha_\gamma < \gamma$. For any $T \in \mathcal{C}_\gamma$, we have $\inf\limits_{\nu \geq 1}\mathbb{E}_{\nu, f_C}[T]\geq \gamma$ and 
    $\inf\limits_{\nu \geq 1}\mathbb{P}_{\nu, f_C}[T \geq \nu] > 0$,
    \begin{align}
        \inf\limits_{\nu \geq 1}\mathbb{P}_{\nu, f_C}[T\leq \nu +\alpha_\gamma|T\geq \nu] \leq \frac{\alpha_\gamma}{\gamma},
    \end{align}
    because otherwise $\inf\limits_{\nu \geq 1}\mathbb{P}_{\nu, f_C}[T \geq \nu +\alpha_\gamma | T \geq \nu] < 1-\nicefrac{\alpha_\gamma}{\gamma}$ 
    with $\inf\limits_{\nu \geq 1}\mathbb{P}_{\nu, f_C}[T\geq \nu] > 0$, implying that $\inf\limits_{\nu \geq 1}\mathbb{E}_{\nu, f_C}[T] \leq \gamma$.

    Let $a = \log\gamma^{(1-\epsilon)}$, then
    \begin{align}
        &e^{a} \mathbb{P}_{\nu, f_C}[\nu \leq T< \nu+\alpha_\gamma | T\geq \nu] \notag\\
        &\leq e^{a} \frac{\alpha_\gamma}{\gamma} = \frac{\alpha_\gamma}{\gamma^\epsilon} \rightarrow 0, \text{ as } \gamma \rightarrow \infty.
    \end{align}

    We then show that the second term at the right-hand side of Eq.~\eqref{eq:lower-proof-break-2} also converges to $0$ as $\gamma \rightarrow \infty$.
    \begin{align}
        &\mathbb{P}_{\nu, f_B}\Big[\max\limits_{\nu\leq j< \nu+\alpha_\gamma}\log\Big(\frac{\mathbb{P}_{\nu, f_B}(H_j)}{\mathbb{P}_{\nu, f_C}(H_j)}\Big) > a\Big|T\geq \nu\Big]\notag\\
        &= \mathbb{P}_{\nu, f_B}\Big[\max\limits_{\nu\leq j< \nu+\alpha_\gamma}\sum_{i=\nu}^j \log\Big(\frac{f_B(X_i)}{f_C(X_i)}\Big) > a |T\geq \nu\Big]\\
        &\overset{(a)}{=}\mathbb{P}_{\nu, f_B}\Big[\max\limits_{\nu\leq j< \nu+\alpha_\gamma}\sum_{i=\nu}^j \log\Big(\frac{f_B(X_i)}{f_C(X_i)}\Big)> a\Big]\\
        &\overset{(b)}{\leq} \mathbb{P}_{\nu, f_B}\Big[\max\limits_{\nu\leq j< \nu+\alpha_\gamma}\sum_{i=\nu}^j \log\Big(\frac{f_B(X_i)}{f_C(X_i)}\Big)\notag\\
        &\qquad\qquad\qquad\qquad> \alpha_\gamma \ (D_{KL}(f_B||f_C)+\epsilon)\Big] \\
        &\rightarrow 0, \text{ as } \gamma \rightarrow \infty.\label{eq:lower-proof-part-2-end}
    \end{align}
    where equality (a) is due to the fact that the event $\{T \geq \nu\}$ is independent from $X_i, \forall i \geq \nu$; inequality (b) is because 
    $a \geq \alpha_\gamma \times (D_{KL}(f_B||f_C)+\epsilon)$;
    and the last step is by applying Lemma~\ref{lemma:WLLN}.

    By Eq.~\eqref{eq:lower-proof-part-1-begin}-\eqref{eq:lower-proof-part-1-end}, we show that 
    \begin{align}
        WADD_{f_B}(T) \geq \frac{\log \gamma}{D_{KL}(f_B||f_0) +o(1)}(1-o(1)),
    \end{align}
    and by Eq.~\eqref{eq:lower-proof-part-2-begin}-\eqref{eq:lower-proof-part-2-end}, we show that 
    \begin{align}
         WADD_{f_B}(T) \geq \frac{\log \gamma}{D_{KL}(f_B||f_C) +o(1)}(1-o(1)).
    \end{align}
    Therefore, we have that
    \begin{align}
        &WADD_{f_B}(T) \notag\\
        &\geq \max\Big\{\frac{\log(\gamma)(1-o(1)) }{D_{KL}(f_B||f_0) +o(1)}, \frac{\log(\gamma)(1-o(1))}{D_{KL}(f_B||f_C) +o(1)}\Big\}\\
        &=\frac{\log \gamma}{\min\{D_{KL}(f_B||f_0), D_{KL}(f_B||f_C)\}+o(1)}(1-o(1)).
    \end{align}

\end{proof}

\section{proof of Theorem~\ref{thm:SCuSum-lower}}\label{sec:proof-SCuSum-lower}

\begin{proof}
    To assist our analysis, we first define intermediate stopping times:
\begin{align}
    &T_{\text{CuSum}_W} := \inf\{t\geq 1: \text{CuSum}_W[t] \geq b_0\},\label{eq:T-CuSumW-def}\\
    &T_{\text{CuSum}^{\texttt{S}}_\Lambda} :=\inf\{t\geq 1: \text{CuSum}^{\texttt{S}}_\Lambda[t] \geq b_C\},\label{eq:T-CuSum-S-Lambda-def}\\
    &T_{\text{CuSum}_{f_B, f_C}} :=\inf\{t\geq 1: \text{CuSum}_{f_B, f_C}[t] \geq b_C\}.\label{eq:T-CuSum-fBfC-def} 
\end{align}
By the algorithmic property of \texttt{S-CuSum}, we have that
\begin{align}
    T_{\texttt{S-CuSum}} &= T_{\text{CuSum}^{\texttt{S}}_\Lambda}\\
    &=T_{\text{CuSum}_W} + T_{\text{CuSum}_{f_B, f_C}}. \label{eq:S-CuSum-stopping-relation}
\end{align}
Hence, if we show 
\begin{align}
    &\mathbb{E}_{\infty}[T_{\texttt{S-CuSum}}] \geq \mathbb{E}_{\infty}[T_{\text{CuSum}_W}] \geq \gamma,\label{eq:SCuSum-false-pre}\\
    &\inf_{\nu\geq 1}\mathbb{E}_{\nu, f_C}[T_{\texttt{S-CuSum}}] \notag\\
    &\geq \min\{\nu + \mathbb{E}_{1, f_C}[T_{\text{CuSum}_{f_B, f_C}}], \mathbb{E}_{\infty}[T_{\text{CuSum}_W}]\} \geq \gamma,\label{eq:SCuSum-false-confusing}
\end{align}
then we have $T_{\texttt{J-CuSum}} \in \mathcal{C}_{\gamma}$.

In the following, we first lower bound the average run time to false alarm for pre-change of $\text{CuSum}_W$ by relating it to a Shiryaev-Robert test. We define the Shiryaev-Roberts statistics corresponding to $\text{CuSum}_W$ as
\begin{align}\label{eq:lower-begin}
    R_W[t] = \sum_{i=1}^t \prod_{j=i}^t \frac{f_B(X_j)}{f_0(X_j)}
\end{align}
with recursion 
\begin{align}
    R_W[t+1] &= (1+R_W[t])\frac{f_B(X_{t+1})}{f_0(X_{t+1})},\\
    R_W[0] &= 0.
\end{align}
Let the stopping time
\begin{align}
    T_{R_W} := \inf\{t\geq 1: R_W[t] \geq e^{b_0}\}.
\end{align}
We have
\begin{align}
    T_{\text{CuSum}_W} \geq  T_{R_W},
\end{align}
because
\begin{align}
    R_W[t] \geq \exp\left\{\max_{1\leq i\leq t}\sum_{j=i}^t W[j]\right\}. 
\end{align}
It then suffices to just lower bound the average run time to false alarm for pre-change of ${R_W}$.

We have that $\{R_W[t] - t\}$ is a martingale w.r.t. $\sigma(X_1,...,X_{t})$ since
\begin{align}
    &\mathbb{E}_{\infty}[(R_W[t+1] - (t+1))|\sigma(X_1,...,X_{t})]\notag\\
    &\overset{(a)}{=} \mathbb{E}_{\infty}\left[(1+R_W[t])\frac{f_B(X_{t+1})}{f_0(X_{t+1})}\Big|\sigma(X_1,...,X_{t})\right]\notag\\
    &\qquad- (t+1)\\
    &\overset{(b)}{=}(1+R_W[t])\mathbb{E}_{\infty}\left[\frac{f_B(X_{t+1})}{f_0(X_{t+1})}\right]-(t+1) \\
    &=(1+R_W[t]) \times 1 -(t+1)\notag\\
    &\overset{(c)}{=}R_W[t] -t, 
\end{align}
where equality (a) is due to the recursion definition of $R_W$, equality (b) is because $R_W[t]$ is $\sigma(X_1,...,X_{t}$-measurable, and equality (c) shows that $\{R_W[t] - t\}$ satisfies the definition of a martingale.

Since $\{R_W[t] - t\}$ is a martingale, by the Doob's optional stopping/sampling theorem, 
\begin{align}
    \mathbb{E}_{\infty}[R_W[T_{R_W}]-T_{R_W}] = \mathbb{E}_{\infty}[R_W[0]] = 0.
\end{align}
Hence, by letting $b_0 = \log \gamma$, we have
\begin{align}
    &\mathbb{E}_\infty[T_{\text{CuSum}_W}] \geq \mathbb{E}_{\infty}[T_{R_W}] \\
    &= \mathbb{E}_{\infty}[R_W[T_{R_W}]] \geq e^{b_0} = \gamma. \label{eq:CuSumW-lower}
\end{align}

In the following, we lower bound the shortest average run time to false alarm for confusing change of $T_{\texttt{S-CuSum}}$. Note that when a false alarm for confusing change is triggered by \texttt{S-CuSum}, there are only two possible cases. In one case, \texttt{S-CuSum} starts updating $T_{\text{CuSum}^{\texttt{S}}_\Lambda}$ after the change point $\nu$, i.e., $\text{CuSum}_W[\nu] < b_0$; in this case, 
\begin{align}
    \mathbb{E}_{\nu, f_C}[T_{\texttt{S-CuSum}}] \geq \nu + \mathbb{E}_{1, f_C}[T_{\text{CuSum}_{f_B, f_C}}].
\end{align}
In the other case, \texttt{S-CuSum} starts updating $T_{\text{CuSum}^{\texttt{S}}_\Lambda}$ after the change point $\nu$, i.e., $\text{CuSum}_W$ passes the threshold before the change point, but $\text{CuSum}^{\texttt{S}}_\Lambda$ passes the threshold after the change point $\nu$; and hence
\begin{align}
    \mathbb{E}_{\nu, f_C}[T_{\texttt{S-CuSum}}] \geq  \mathbb{E}_{\infty}[T_{\text{CuSum}_W}].
\end{align}
Therefore, the shortest average run time to false alarm for confusing change $\inf_{\nu\geq 1}\mathbb{E}_{\nu, f_C}[T_{\texttt{S-CuSum}}]$ is lower bounded by $\min\{\nu + \mathbb{E}_{1, f_C}[T_{\text{CuSum}_{f_B, f_C}}], \mathbb{E}_{\infty}[T_{\text{CuSum}_W}]\}$.

In the following, we lower bound $\mathbb{E}_{1, f_C}[T_{\text{CuSum}_{f_B, f_C}}]$ following the similar reasoning as for lower bounding $\mathbb{E}_{\infty}[T_{\text{CuSum}_W}]$ in Eq.~\eqref{eq:lower-begin}-\eqref{eq:CuSumW-lower}. Specifically, we let 
\begin{align}\label{eq:SCuSum-lower-part2-begin}
    &R_{\Lambda} := \sum_{i=1}^t\prod_{j=i}^t \frac{f_B(X_j)}{f_C(X_j)},\\
    &T_{R_\Lambda} :=\inf\left\{t \geq 1: R_{\Lambda}[t] \geq e^{b_C}\right\},
\end{align}
and we have that $\{R_{\Lambda}[t] - t\}$ is a martingale w.r.t. $\sigma(X_1,...,X_{t})$.
Then by Doob's optional stopping theorem, we have 
\begin{align}
    &\mathbb{E}_{1, f_C}[T_{\text{CuSum}_{f_B, f_C}}] \geq \mathbb{E}_{1, f_C}[T_{R_\Lambda}] \\
    &= \mathbb{E}_{1, f_C}[R_\Lambda[T_{R_\Lambda}]] \geq e^{b_C} = \gamma, \label{eq:CuSumLambda-lower}
\end{align}
where the last inequality is by letting $b_C = \log \gamma$.

By Eq.~\eqref{eq:SCuSum-false-pre}\eqref{eq:SCuSum-false-confusing}\eqref{eq:CuSumW-lower}\eqref{eq:CuSumLambda-lower}, we show that $T_{\texttt{J-CuSum}} \in \mathcal{C}_{\gamma}$.

\end{proof}

\section{proof of Theorem~\ref{thm:SCuSum-upper}}\label{sec:proof-SCuSum-upper}

\begin{proof}
    By the fact that $\text{CuSum}_W[t]$ and $\text{CuSum}^{\texttt{S}}_\Lambda[t]$ are always non-negative and are zero when $t=0$, the worse-case of the average detection delay happens when change point $\nu = 1$. 

    To assist the analysis, we will use intermediate stopping times $T_{\text{CuSum}_W}$, $T_{\text{CuSum}^{\texttt{S}}_\Lambda}$, and $T_{\text{CuSum}_{f_B, f_C}}$ defined in Eq.~\eqref{eq:T-CuSumW-def}\eqref{eq:T-CuSum-S-Lambda-def}\eqref{eq:T-CuSum-fBfC-def} respectively. 

    By the algorithmic properties of \texttt{S-CuSum}, we have that
    \begin{align}
        \mathbb{E}_{1, f_B}[T_{\texttt{S-CuSum}}] = \mathbb{E}_{1, f_B}[T_{\text{CuSum}_W}] + \mathbb{E}_{1, f_B}[T_{\text{CuSum}_{f_B, f_C}}].
    \end{align}

    Let $\alpha_{b_0} = \nicefrac{b_0}{D_{KL}(f_B||f_0)}$ and $\alpha_{b_C} = \nicefrac{b_C}{D_{KL}(f_B||f_C)}$. Then,
    \begin{align}
        &\mathbb{E}_{1, f_B}[\nicefrac{T_{\text{CuSum}_W}}{\alpha_{b_0}}]\notag\\
        &=\sum_{\ell=0}^\infty \mathbb{P}_{1, f_B}[\nicefrac{T_{\text{CuSum}_W}}{\alpha_{b_0}} > \ell]\\
        &\leq 1 + \sum_{\ell=1}^\infty \mathbb{P}_{1, f_B}[\nicefrac{T_{\text{CuSum}_W}}{\alpha_{b_0}} > \ell]\\
        &= 1+ \sum_{\ell=1}^\infty \mathbb{P}_{1, f_B}[\forall 1\leq t\leq \ell\alpha_{b_0}: \text{CuSum}_W[t] < b_0]\\
        &\leq 1+ \sum_{\ell=1}^\infty \bigcap_{k=1}^\ell \mathbb{P}_{1, f_B}[\text{CuSum}_W[k\alpha_{b_0}] < b_0]\\
        &\overset{(a)}{\leq} 1+ \sum_{\ell=1}^\infty \prod_{k=1}^\ell \mathbb{P}_{1, f_B}\left[\max_{i: (k-1)\alpha_{b_0}+1 \leq i \leq k \alpha_{b_0}} \sum_{j=i}^{k \alpha_{b_0}} \frac{f_B(X_j)}{f_0(X_j)}< b_0\right],\label{eq:upper-intermideate-step-1}
    \end{align}
    \begin{align}
        &\mathbb{E}_{1, f_B}[\nicefrac{T_{\text{CuSum}_{f_B, f_C}}}{\alpha_{b_C}}]\notag\\
        &=\sum_{\ell=0}^\infty\mathbb{P}_{1, f_B}[\nicefrac{T_{\text{CuSum}_{f_B, f_C}}}{\alpha_{b_C}} > \ell]\\
        &\leq 1+\sum_{\ell=1}^\infty\mathbb{P}_{1, f_B}[\nicefrac{T_{\text{CuSum}_{f_B, f_C}}}{\alpha_{b_C}} > \ell]\\
        &= 1 +\sum_{\ell=1}^\infty \mathbb{P}_{1, f_B}[\forall 1\leq t\leq \ell\alpha_{b_C}: \text{CuSum}_{f_B, f_C}[t] < b_C]\\
        &\leq 1 +\sum_{\ell=1}^\infty\bigcap_{k=1}^\ell \mathbb{P}_{1, f_B}[\text{CuSum}_{f_B, f_C}[k\alpha_{b_C}] < b_C]\\
        &\overset{(a)}{\leq} 1 +\sum_{\ell=1}^\infty\prod_{k=1}^\ell \mathbb{P}_{1, f_B}\left[\max_{i: (k-1)\alpha_{b_0}+1 \leq i \leq k \alpha_{b_0}} \sum_{j=i}^{k \alpha_{b_0}} \frac{f_B(X_j)}{f_C(X_j)} < b_C\right] \label{eq:upper-intermideate-step-2}
    \end{align}
    where inequalities (a) are by the definitions of $\text{CuSum}_W$ and $\text{CuSum}_{f_B, f_C}$ and by the independency among the random variables.  

    It follows from Lemma~\ref{lemma:WLLN} that 
    \begin{align}
        &\frac{\max\limits_{i: 1 \leq i \leq \alpha_{b_0}} \sum_{j=i}^{ \alpha_{b_0}} \frac{f_B(X_j)}{f_0(X_j)}}{b_0} \overset{p}{\rightarrow} \beta,\\
        &\frac{\max\limits_{i: 1 \leq i \leq \alpha_{b_C}} \sum_{j=i}^{ \alpha_{b_0}} \frac{f_B(X_j)}{f_C(X_j)}}{b_C} \overset{p}{\rightarrow} \beta,
    \end{align}
    where $\beta > 1$. Therefore, as $b_0, b_C\rightarrow \infty$,
    \begin{align}
        &\mathbb{P}_{1, f_B}\left[\max\limits_{i: 1 \leq i \leq \alpha_{b_0}} \sum_{j=i}^{ \alpha_{b_0}} \frac{f_B(X_j)}{f_0(X_j)} < b_0\right] \rightarrow 0,\\
        &\mathbb{P}_{1, f_B}\left[\max\limits_{i: 1 \leq i \leq \alpha_{b_0}} \sum_{j=i}^{ \alpha_{b_0}} \frac{f_B(X_j)}{f_C(X_j)} < b_C\right] \rightarrow 0.
    \end{align}
    This implies that 
    \begin{align}
        &\mathbb{P}_{1, f_B}\left[\max\limits_{i: 1 \leq i \leq \alpha_{b_0}} \sum_{j=i}^{ \alpha_{b_0}} \frac{f_B(X_j)}{f_0(X_j)} < b_0\right] \leq \delta, \label{eq:prob-small-1}\\
        &\mathbb{P}_{1, f_B}\left[\max\limits_{i: 1 \leq i \leq \alpha_{b_0}} \sum_{j=i}^{ \alpha_{b_0}} \frac{f_B(X_j)}{f_C(X_j)} < b_C\right] \leq \delta, \label{eq:prob-small-2}
    \end{align}
    where $\delta$ can be arbitrarily small for large $b_0$ and $b_C$.

    By Eq.~\eqref{eq:prob-small-1}\eqref{eq:prob-small-2} and Eq.~\eqref{eq:upper-intermideate-step-1}\eqref{eq:upper-intermideate-step-2}, we have
    \begin{align}
        &\mathbb{E}_{1, f_B}[\nicefrac{T_{\text{CuSum}_W}}{\alpha_{b_0}}] + \notag\\
        &\leq 1 + \sum_{\ell=1}^\infty \delta^\ell\\
        &=\frac{1}{1-\delta},
    \end{align}
    \begin{align}
        &\mathbb{E}_{1, f_B}[\nicefrac{T_{\text{CuSum}_{f_B, f_C}}}{\alpha_{b_C}}]\notag\\
        &\leq 1 + \sum_{\ell=1}^\infty \delta^\ell\\
        &=\frac{1}{1-\delta}.
    \end{align}
    This implies that, as $\gamma \rightarrow \infty$,
    \begin{align}
        &\mathbb{E}_{1, f_B}[T_{\texttt{S-CuSum}}] \notag\\
        &\leq \frac{\alpha_{b_0} + \alpha_{b_C}}{1-\delta}\\
        & \leq \left(\frac{\log \gamma}{D_{KL}(f_B||f_0)}+\frac{\log\gamma}{D_{KL}(f_B||f_C)}\right)(1+o(1)), \notag\\
        &\qquad\qquad\qquad\qquad\qquad\qquad\text{ if } b_0 = b_C = \log \gamma\\
        & \leq \frac{2\log\gamma}{\min\{D_{KL}(f_B||f_0), D_{KL}(f_B||f_C)\}}(1+o(1)), \notag\\
        &\qquad\quad\text{ if } b_0 = \frac{D_{KL}(f_B||f_0)}{\min\{D_{KL}(f_B||f_0), D_{KL}(f_B||f_C)\}},\notag\\
        &\qquad\qquad b_C = \frac{D_{KL}(f_B||f_C)}{\min\{D_{KL}(f_B||f_0), D_{KL}(f_B||f_C)\}}.
    \end{align}
\end{proof}

\section{proof of Theorem~\ref{thm:JCuSum-lower}}\label{sec:proof-JCuSum-lower}

\begin{proof}
To assist the analysis, we will use intermediate stopping times $T_{\text{CuSum}_W}$ and $T_{\text{CuSum}_{f_B, f_C}}$ defined in Eq.~\eqref{eq:T-CuSumW-def}\eqref{eq:T-CuSum-fBfC-def} respectively and 
    \begin{align}\label{eq:JCuSum_lambda_def}
        T_{\text{CuSum}_{\Lambda}^{\texttt{J}}} := \inf\{t \geq 1: \text{CuSum}_{\Lambda}^{\texttt{J}}[t] \geq b_C\}.
    \end{align}

By algorithmic property of \texttt{J-CuSum}, we have that
\begin{align}
    T_{\texttt{J-CuSum}} = \max\{T_{\text{CuSum}_W}, T_{\text{CuSum}_\Lambda^{\texttt{J}}}\}. 
\end{align}
Hence, if we show 
\begin{align}
    &\mathbb{E}_{\infty}[T_{\texttt{J-CuSum}}] \geq \mathbb{E}_{\infty}[T_{\text{CuSum}_W}] \geq \gamma,\label{eq:JCuSum-false-pre}\\
    &\inf_{\nu\geq 1}\mathbb{E}_{\nu, f_C}[T_{\texttt{J-CuSum}}] \geq \inf_{\nu\geq 1}\mathbb{E}_{\nu, f_C}[T_{\text{CuSum}_\Lambda^{\texttt{J}}}] \geq \gamma,\label{eq:JCuSum-false-confusing}
\end{align}
then we have $T_{\texttt{J-CuSum}} \in \mathcal{C}_{\gamma}$.

By the analysis in Appendix~\ref{sec:proof-SCuSum-lower}, i.e., Eq.~\eqref{eq:lower-begin}-\eqref{eq:CuSumW-lower}, we already have that  
\begin{align}
    \mathbb{E}_{\infty}[T_{\text{CuSum}_W}] \geq \gamma.
\end{align}

In the following, we lower bound the shortest average run time to false alarm for confusing change of $T_{\texttt{J-CuSum}}$. To assist this analysis, we define $T_{\text{CuSum}_W}$ in terms of a sequence of sequential probability ratio tests as in~\cite{siegmund1985sequential}. Specifically, let 
    \begin{align}
        S[t] &:= \sum_{i=1}^t W[i]\label{eq:siegmund-begin}\\
        N_1 &:= \inf\{t\geq 1: S[t] \notin (0, b_0)\},\\
        N_2 &:= \inf\{t\geq 1: S[N_1+t] -S[N_1] \notin (0, b_0)\},\\
        N_k &:=\inf\{t\geq 1: S[N_1+N_2+...+N_{k-1}+t]\notag\\
        &\qquad\quad-S[N_1+N_2+...+N_{k-1}] \notin (0, b_0)\},\\
        M &:=\inf\{k\geq 1: S[N_1+N_2+...+N_k]\notag\\
        &\qquad\quad-S[N_1+N_2+...+N_{k-1}] \geq b_0\}; 
    \end{align}
this way, 
\begin{align}
    T_{\text{CuSum}_W} \equiv N_1 + N_2 + ... + N_M;
\end{align}
besides, we also let
    \begin{align}
        \alpha &:= \mathbb{P}_{f_0}[S[N_1]\geq b_0],\\
        \beta &:= \mathbb{P}_{f_B}[S[N_1]\leq a], a \rightarrow 0^{-}.\label{eq:siegmund-end}
    \end{align}

By the algorithmic property of $\texttt{J-CuSum}$, when confusing change occurs and $\nu \rightarrow \mathbb{E}_{f_0}[N_1]^{-}$, $\texttt{J-CuSum}$ has the shortest average run time for $\text{CuSum}_\Lambda^{\texttt{J}}$ to pass the threshold; and
\begin{align}
    &\mathbb{E}_{f_0}[\text{CuSum}_\Lambda^{\texttt{J}}[\nu]]\notag\\
    &\overset{(a)}{=}\mathbb{E}_{f_0}\Big[\log\frac{f_B(X)}{f_C(X)}\Big]\times \nu\\
    &\simeq \mathbb{E}_{f_0}\Big[\log\frac{f_B(X)}{f_C(X)}\Big]\times\mathbb{E}_{f_0}[N_1], \text{ as } \nu \rightarrow \mathbb{E}_{f_0}[N_1]^{-}\\
    &\overset{(b)}{\simeq} \mathbb{E}_{f_0}\Big[\log\frac{f_B(X)}{f_C(X)}\Big]\times\frac{a(e^{b_0}-1)+b_0(1-e^a)}{-(e^{b_0}-e^{a})D_{KL}(f_0||f_B)}, a\rightarrow 0^{-}\\
    &\simeq 0 \label{eq:JCuSum-approx}
\end{align}
where equality (a) is by Wald's identity; approximation (b) is by~\cite[Eq. (2.15)]{siegmund1985sequential}. 
Hence  
\begin{align}
    &\inf\limits_{\nu \geq 1} \mathbb{E}_{\nu, f_C}[T_{\text{CuSum}_\Lambda^{\texttt{J}}}]\notag\\
    &=\inf\limits_{\nu \geq 1} \sum_{t=0}^\infty \mathbb{P}_{\nu, f_C}[T_{\text{CuSum}_\Lambda^{\texttt{J}}} \geq t]\\
    &=\inf\limits_{\nu \geq 1} \sum_{t=0}^\infty \mathbb{P}_{\nu, f_C}[{\text{CuSum}_\Lambda^{\texttt{J}}[t]} < b_C]\\
    &\overset{(a)}{\simeq}  \sum_{t=0}^\infty \mathbb{P}_{f_C}[\text{CuSum}_{f_B, f_C}[t] < b_C]\\
    &=\sum_{t=0}^\infty \mathbb{P}_{f_C}[T_{\text{CuSum}_{f_B, f_C}} \geq t]\\
    &=\mathbb{E}_{f_C}[T_{\text{CuSum}_{f_B, f_C}}]\\
    &\overset{(b)}{\geq} \gamma,
\end{align}
where approximation (a) follows from Eq.~\eqref{eq:JCuSum-approx}, and inequality (b) follows from the analysis in Appendix~\ref{sec:proof-SCuSum-lower}, i.e., Eq.~\eqref{eq:SCuSum-lower-part2-begin}-\eqref{eq:CuSumLambda-lower}.

\end{proof}

\begin{figure}[t]
    \begin{minipage}[b]{0.48\linewidth}
        \includegraphics[width=\columnwidth]{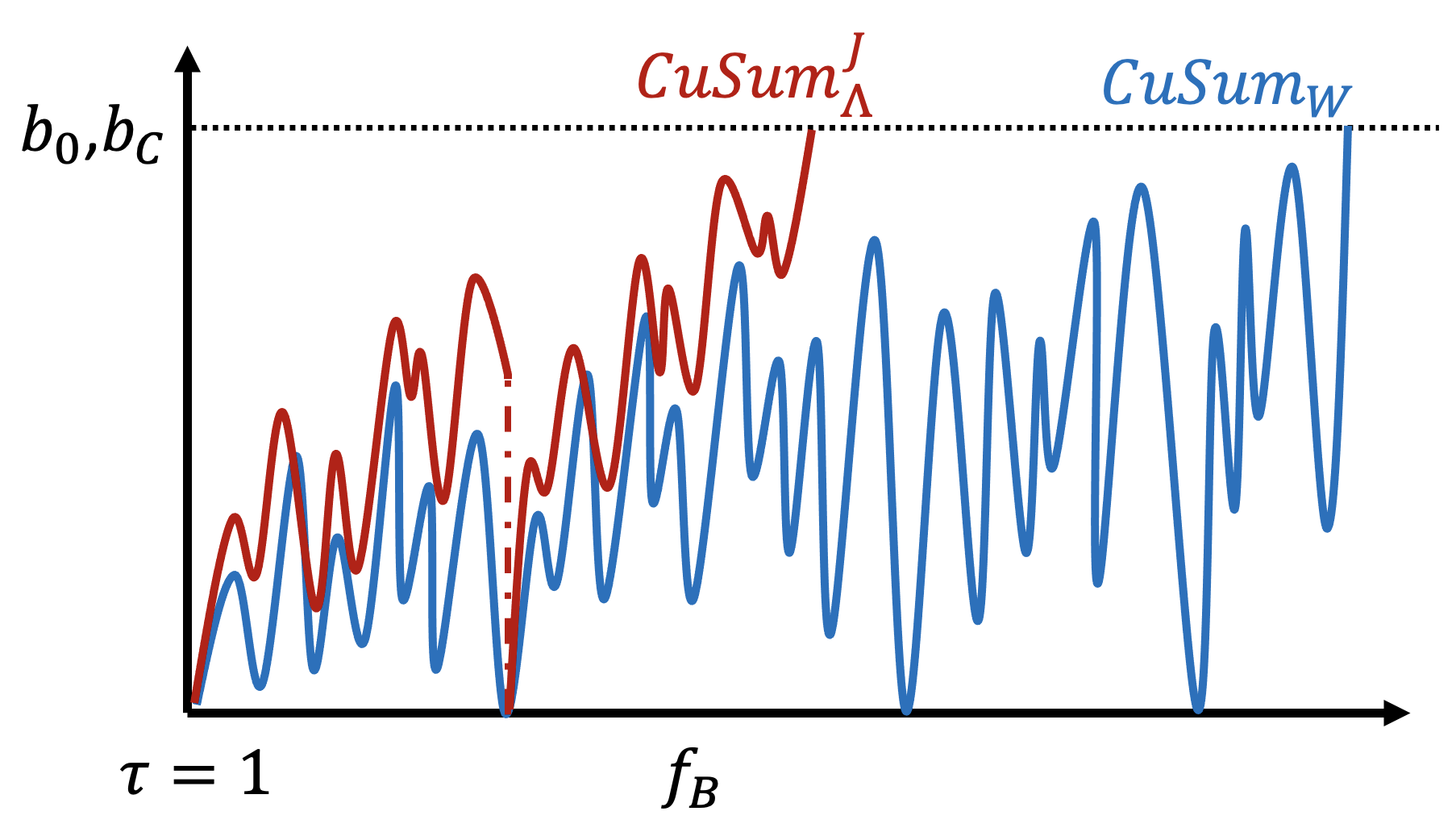} 
        \subcaption{Case 1}
        \label{fig:detection-delay-scenario-1}
    \end{minipage}
    \hfill
    \begin{minipage}[b]{0.48\linewidth}
        \includegraphics[width=\columnwidth]{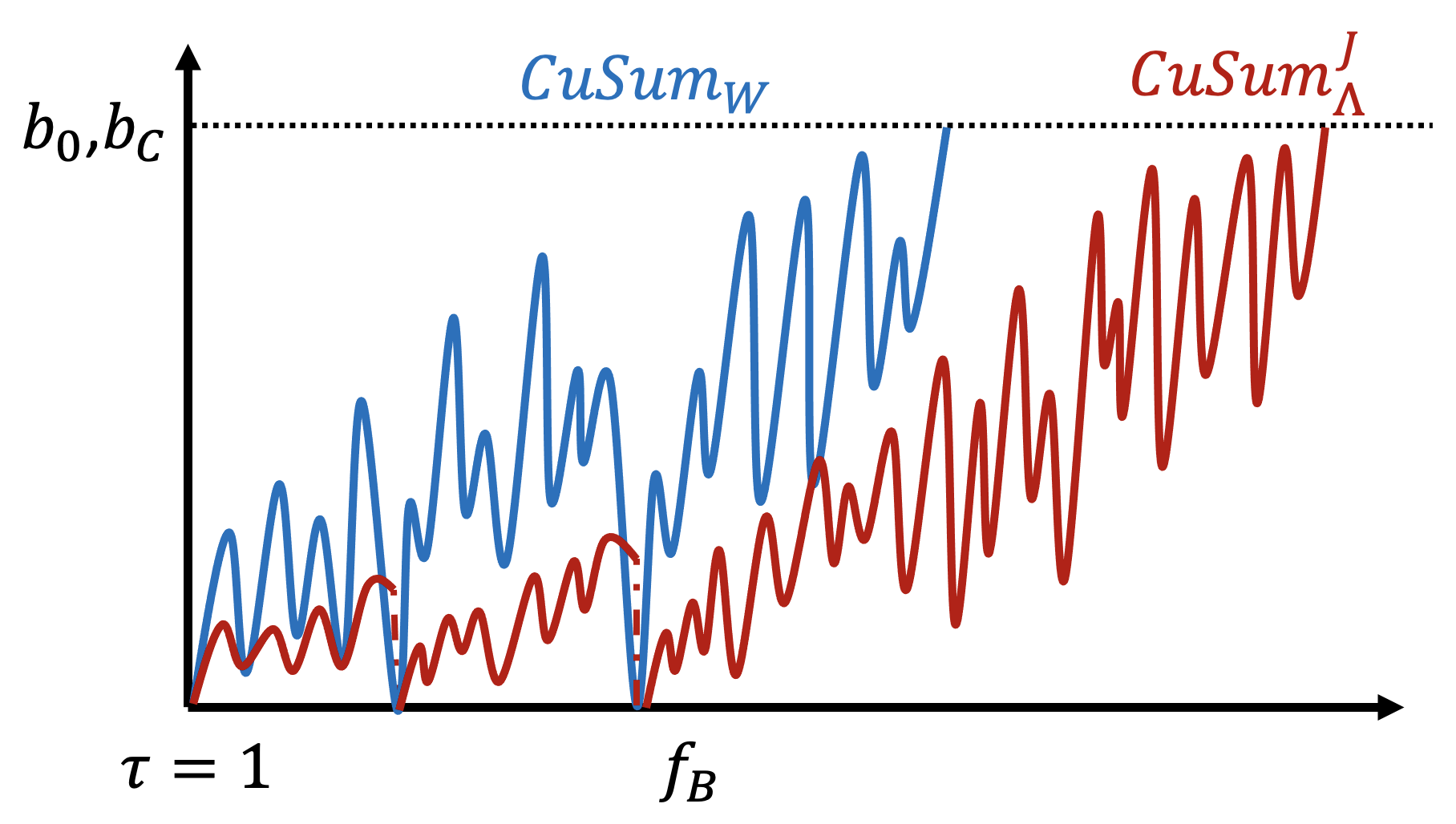}
        \subcaption{Case 2}
        \label{fig:detection-delay-scenario-2}
    \end{minipage}
    \caption{Illustration of Possible Performances of \texttt{J-CuSum} when $\tau =1$: in Case 1, $\text{CuSum}^{\texttt{J}}_\Lambda$ passes its threshold \textit{before} $\text{CuSum}_W$ does. In Case 2, $\text{CuSum}^{\texttt{J}}_\Lambda$ passes its threshold \textit{after} $\text{CuSum}_W$ does.}
\end{figure}

\section{proof of Theorem~\ref{thm:JCuSum-upper}}\label{sec:proof-JCuSum-upper}

\begin{proof}
     By the fact that $\text{CuSum}_W[t]$ and $\text{CuSum}^{\texttt{J}}_\Lambda[t]$ are always non-negative and are zero when $t=0$, the worse-case of the average detection delay happens when change point $\nu = 1$. 

    To assist the analysis, we will use intermediate stopping times $T_{\text{CuSum}_W}$ and $T_{\text{CuSum}_{f_B, f_C}}$ defined in Eq.~\eqref{eq:T-CuSumW-def}\eqref{eq:T-CuSum-fBfC-def} respectively and $T_{\text{CuSum}_{\Lambda}^{\texttt{J}}}$ defined in Eq.~\eqref{eq:JCuSum_lambda_def}.

    By the algorithmic property of \texttt{J-CuSum}, we have that,
    when $\text{CuSum}^{\texttt{J}}_\Lambda$ passes threshold $b_C$, there are only two possible cases: 1) $\text{CuSum}_W$ has not passed threshold $b_0$ yet (as illustrated in Figure~\ref{fig:detection-delay-scenario-1}), 2) $\text{CuSum}_W$ has already passed threshold $b_0$ (as illustrated in Figure~\ref{fig:detection-delay-scenario-2}). In the following, we analyze each case separately.

    We first consider case 1, the case that $\text{CuSum}_\Lambda$ passes threshold $b_C$ while $\text{CuSum}_W$ has not passed threshold $b_0$ (as illustrated in Figure~\ref{fig:detection-delay-scenario-1}). In this case, the worse case detection delay of \texttt{J-CuSum} is simply upper bounded by that of $\text{CuSum}_W$. And by the analysis in Appendix~\ref{sec:proof-SCuSum-upper}, we have that, with $b_0 = \log \gamma$, $\gamma \rightarrow \infty$, in case 1:
    \begin{align}
        &\mathbb{E}_{1, f_B}[T_{\texttt{J-CuSum}}] \notag\\
        &\leq \mathbb{E}_{1, f_B}[T_{\text{CuSum}_W}]\\
        &\leq \frac{\log \gamma}{D_{KL}(f_B||f_0)}(1+o(1)).
    \end{align}

    We then proceed to consider case 2, the case that $\text{CuSum}_\Lambda$ passes threshold $b_C$ when $\text{CuSum}_W$ has already passed threshold $b_0$ (as illustrated in Figure~\ref{fig:detection-delay-scenario-2}). To assist this analysis, as introduced in Eq.~\eqref{eq:siegmund-begin}-\eqref{eq:siegmund-end}, we follow~\cite{siegmund1985sequential} defining $T_{\text{CuSum}_W}$ in terms of a sequence of sequential probability ratio tests.
    Then we have that, with $b_0 = b_C =  \log \gamma$, $\gamma \rightarrow \infty$, in case 2: 
    \begin{align}
        &\mathbb{E}_{1, f_B}[T_{\texttt{J-CuSum}}] \notag\\
        &\leq \mathbb{E}_{1, f_B}[T_{\text{CuSum}^{\texttt{J}}_\Lambda}]\\
        &\leq \mathbb{E}_{1, f_B}\left[\sum_{k=1}^{M-1}N_k\right] + \mathbb{E}_{1, f_B}[T_{\text{CuSum}_{f_B, f_C}}]\\
        &= \mathbb{E}_{1, f_B}\left[(M-1)\cdot N_1\right] + \mathbb{E}_{1, f_B}[T_{\text{CuSum}_{f_B, f_C}}]\\
        &\overset{(a)}{=}\left(\frac{1}{\mathbb{P}_{1, f_B}[S[N_1]\geq b_0]}-1\right) \mathbb{E}_{1, f_B}\left[N_1\right]\notag\\
        &\qquad+ \mathbb{E}_{1, f_B}[T_{\text{CuSum}_{f_B, f_C}}]\\
        &\overset{(b)}{\simeq} \left(\frac{1}{1-\beta}-1 \right)\mathbb{E}_{1, f_B}\left[N_1\right]+ \mathbb{E}_{1, f_B}[T_{\text{CuSum}_{f_B, f_C}}]\\
        &\overset{(c)}{\simeq} \left(\frac{e^a(e^{b_0}-1)}{e^{b_0}(1-e^a)}\right)\mathbb{E}_{1, f_B}\left[N_1\right]\notag\\
        &\qquad+ \mathbb{E}_{1, f_B}[T_{\text{CuSum}_{f_B, f_C}}], a\rightarrow 0^{-}\\
        &\overset{(d)}{\simeq} \frac{e^a(e^{b_0}-1)[ae^a(e^{b_0}-1)+b_0e^{b_0}(1-e^a)]}{e^{b_0}(1-e^a)(e^{b_0}-e^a)D_{KL}(f_B||f_0)}\notag\\
        &\qquad+ \mathbb{E}_{1, f_B}[T_{\text{CuSum}_{f_B, f_C}}], a\rightarrow 0^{-}\label{eq:plug-EN}\\
        &\overset{(e)}{=} \frac{b_0 -1 +\nicefrac{1}{e^{b_0}}}{D_{KL}(f_B||f_0)}+ \mathbb{E}_{1, f_B}[T_{\text{CuSum}_{f_B, f_C}}]\\
        &\overset{(f)}{\leq} \left(\frac{\log \gamma}{D_{KL}(f_B||f_0)}+\frac{\log \gamma}{D_{KL}(f_B||f_C)}\right)(1+o(1)).
    \end{align} 
    where equality (a) is by Wald's identity, approximation (b) is following~\cite[Eq. (2.52)(2.53)]{siegmund1985sequential}, approximation (c) is following~\cite[Eq. (2.11)(2.12)]{siegmund1985sequential}, approximation (d) is by~\cite[Eq. (2.15)]{siegmund1985sequential}, equality (e) is by L'Hôpital's rule, and inequality (f) is by the analysis in Appendix~\ref{sec:proof-SCuSum-upper}.
\end{proof}

\end{document}